\documentclass{amsart}
\usepackage{amsfonts, amsbsy, amsmath, amssymb}
\newtheorem{thm}{Theorem}[section]
\newtheorem{lem}[thm]{Lemma}
\newtheorem{cor}[thm]{Corollary}
\newtheorem{prop}[thm]{Proposition}
\newtheorem{exmp}[thm]{Example}

\newtheorem{algo}[thm]{Algorithm}
\numberwithin{equation}{section}

\theoremstyle{definition}

\newcounter{step}

\renewcommand{\thestep}{\sc Step \arabic{step}}

\begin{document}

\title{On the Number of Inequivalent Binary Self-Orthogonal Codes}

\author{Xiang-dong Hou}
\address{Department of Mathematics,
University of South Florida, Tampa, FL 33620}
\email{xhou@tarski.math.usf.edu}


 
\begin{abstract}
Let $\Psi_{k,n}$ denote the number of inequivalent binary self-orthogonal
$[n,k]$ codes. We present a method which allows us to compute $\Psi_{k,n}$
explicitly for a moderate $k$ and an arbitrary $n$. Included in this paper
are explicit formulas for $\Psi_{k,n}$ with $k\le 5$.
\end{abstract}

\maketitle


\section{Introduction}
Let $\Bbb F_2$ be the binary field. Throughout the paper, a {\em code} is a
subspace of $\Bbb F_2^n$ for some integer $n>0$; an $[n,k]$ code is a
$k$-dimensional subspace of $\Bbb F_2^n$. Let $\langle\cdot,\cdot\rangle$ denote the usual inner product of $\Bbb F_2^n$. A code $C\subset\Bbb F_2^n$
is called self-orthogonal if $C\subset C^\bot$ where $C^\bot=\{x\in\Bbb F_2^n:\langle x,y\rangle=0\ \text{for all}\ y\in C\}$. The number of
$[n,k]$ self-orthogonal codes is known. Let $\Phi_{n,k}$ be the set of all
$[n,k]$ self-orthogonal codes. Then
\[
|\Phi_{n,k}|=
\begin{cases}
\displaystyle
\frac{\prod_{j=1}^k(2^{n+1-2j}-1)}{\prod_{j=1}^k(2^j-1)}&\text{if $n$ is odd},\vspace{4pt}\cr
\displaystyle
\frac{(2^{n-k}-1)\prod_{j=1}^{k-1}(2^{n-2j}-1)}{\prod_{j=1}^k(2^j-1)}&\text{if $n\ge 2$ is even},
\end{cases}
\]
see MacWilliams and Sloane \cite[Ch. 19, \S6]{Mac81}. 

Let $\frak S_n$ be the symmetric group on $\{1,\dots, n\}$. $\frak S_n$ acts on $\Bbb F_2^n$ by permuting the coordinates of $\Bbb F_2^n$. Two codes $C_1, C_2\subset\Bbb F_2^n$ are called equivalent if there is a $\sigma\in\frak S_n$ such that $C_2=C_1^\sigma$. Self-orthogonality
of codes is preserved under the equivalence. Let $\Psi_{k,n}$ be the number
of inequivalent $[n,k]$ self-orthogonal codes. Unlike $|\Phi_{n,k}|$, the number $\Psi_{k,n}$ is much more difficult to compute. In \cite{Con80}, Conway and Pless determined $\Psi_{k,2k}$ for $k\le 15$ after classifying
doubly even $[32,16]$ codes. (Also see \cite{Con92} for an update on $\Psi_{15,30}$.) 
In a recent paper \cite{Hou1}, the author considered the asymptotic behavior of
$\Psi_{k,2k}$ and proved that
\[
\Psi_{k,2k}\sim \tau \cdot(2k)!2^{\frac 12k(k-1)}\quad\text{as}\ k\to\infty,
\]
where $\tau=\prod_{j=1}^\infty(1+2^{-j})$. However, not much else is known about the number 
$\Psi_{k,n}$.

The situation described above brings up a basic question: can $\Psi_{k,n}$ be computed explicitly?
In this paper, we will see that the answer is ``yes'' for a moderate $k$ and an arbitrary $n$.

Our interest in the number $\Psi_{k,n}$ is motivated by the important role it can play
in the classification of $[n,k]$ self-orthogonal codes. Without knowing the number $\Psi_{k,n}$
beforehand, the known method to classify $[n,k]$ self-orthogonal codes relies on the mass
formula. The algorithm of this method is sketched as follows. Assume that a list of pairwise
inequivalent $[n,k]$ self-orthogonal codes $C_1,\dots,C_t$ has been found and the cardinality of the automorphism group of each $C_i$ has been determined. Then
\begin{equation}\label{1.1}
n!\sum_{i=1}^t\frac 1{|{\rm Aut}(C_i)|}
\end{equation}
is the number of $[n,k]$ self-orthogonal codes equivalent to one of $C_1,\dots, C_t$.
If the sum in \eqref{1.1} is $<|\Phi_{n,k}|$, search for an $[n,k]$ self-orthogonal code
$C_{t+1}$ which is inequivalent to all of $C_1,\dots,C_t$, and compute $|{\rm Aut}(C_{t+1})|$.
Add $C_{t+1}$ to the list $C_1,\dots, C_t$ and update the sum in \eqref{1.1}. The list $C_1,\dots, C_t$ is complete when the sum in \eqref{1.1} equals $|\Phi_{n,k}|$. In fact, Conway, Pless and Sloane's classifications of self-orthogonal codes of length up to 30 and doubly even
$[32,16]$ codes were obtained using this method \cite{Con80, Ple72, Ple75}. On the other hand, if the number $\Psi_{k,n}$ is known beforehand, the algorithm to classify $[n,k]$ 
self-orthogonal codes is greatly simplified. One only has to find $\Psi_{k,n}$ pairwise
inequivalent $[n,k]$ self-orthogonal codes.

A key step in the computation of $\Psi_{k,n}$ is to determine the numbers of zeros of certain quadratic forms defined on $\Bbb F_2^n$. This 
requires us to be able to tell the canonical forms of those quadratic forms. For this purpose, a brief review of canonical forms of binary quadratic forms is given in Section 2. 
In Section 3, we outline the method for computing $\Psi_{k,n}$. We also prove a few preliminary results to be used later. We
derive the formula for $\Psi_{3,n}$ in Section 4 and the formula for $\Psi_{4,n}$ in Section 5.
Most of the details of the computations in Sections 4 and 5 are included. In Section 6, we give the formula for 
$\Psi_{5,n}$ but omit the details of the computations. It should be clear from the paper that the
method works for $k$ beyond the range considered here. In Section 7, we give the numerical values of
$\Psi_{k,n}$ for $k\le 5$ and $n\le 40$.

In our notation, $\Bbb N=\{0,1,2,\dots\}$. $M_{k\times n}$ is the set of all $k\times n$ matrices over $\Bbb F_2$. The $n\times n$ identity matrix is denoted by $I_n$ or simply $I$ when $n$ is clear from the context. ${\bf 1}(n)$
is the $1\times n$ all one vector. For two matrices 
$A,B$,
\[
A\oplus B:=\left[
\begin{matrix}
A\cr
&B
\end{matrix}\right].
\]
For each function $f:\Bbb F_2^n\to\Bbb F_2$, $Z(f):=\{x\in\Bbb F_2^n:f(x)=0\}$.  
For $n\in \Bbb Z$, $\nu(n)$ is the 2-adic order of $n$. We also define
\[
\delta(n)=
\begin{cases}
0&\text{if}\ n\le 0,\cr
1&\text{if}\ n> 0.
\end{cases}
\]
A congruence $a\equiv b\pmod n$ is abbreviated as $a\equiv b\ (n)$.


\section{Binary Quadratic Forms}

A quadratic form in $n$ variables over $\Bbb F_2$ is a function $f:\Bbb F_2^n\to\Bbb F_2$ of the form
\begin{equation}\label{3.-1}
f(x_1,\dots,x_n)=\sum_{i\le j}a_{ij}x_ix_j,
\end{equation}
where $a_{ij}\in\Bbb F_2$. Of course, the quadratic form in \eqref{3.-1}
can also be written as
\[
f(x_1,\dots,x_n)=\sum_{i< j}a_{ij}x_ix_j+\sum_i a_{ii}x_i.
\]
Let
\[
\Lambda_n=\{A\in M_{n\times n}: A^T=A\ \text{and all diagonal entries
of $A$ are 0}\}.
\]
Then there is a bijection between $M_{n\times n}/\Lambda_n$ and the set of
all quadratic forms in $n$ variables over $\Bbb F_2$:
\[
A+\Lambda_n\longleftrightarrow (x_1,\dots,x_n)A(x_1,\dots,x_n)^T.
\]
We say that $A$ is a matrix of the quadratic form $(x_1,\dots,x_n)A(x_1,\dots,x_n)^T$.

Let $f(x_1,\dots,x_n)$ and $g(x_1,\dots,x_n)$ be two quadratic forms over
$\Bbb F_2$ with matrices $A$ and $B$ respectively. If $f(x_1,\dots,x_n)=g\bigl((x_1,\dots,x_n)Q\bigr)$ for some $Q\in{\rm GL}(n,
\Bbb F_2)$, we say that $f$ and $g$ are linearly equivalent and we write $f\cong g$. Clearly, $f\cong g$ if and only if
\begin{equation}\label{3.-2}
A\equiv QBQ^T \pmod {\Lambda_n}
\end{equation}
for some $Q\in{\rm GL}(n,\Bbb F_2)$. If two matrices $A,B\in M_{n\times n}$
satisfy \eqref{3.-2}, we write $A\cong B$. Therefore, finding the canonical
form of a quadratic form under linear equivalence is the same as finding
the canonical form of its matrix under the equivalence $\cong$.
In practice, it is more convenient to work with the matrices than the 
quadratic forms themselves.

In the following theorem, we collect some well known results on the canonical forms and numbers of zeros of binary quadratic forms.

\begin{thm}\label{T3.1}\

\begin{enumerate}
\item[(i)]
Every $A\in M_{n\times n}$ is $\cong$ equivalent to exactly one of the following canonical forms:
\[
\left[
\kern-0.8cm
\begin{matrix}
{\scriptstyle 2r}\left\{
\phantom{\begin{matrix} 0\cr 0\cr \vdots\cr 0\cr 0\cr \end{matrix}}\right.
\begin{matrix} 
0&1\cr 
0&0\cr 
&&\ddots\cr 
&&&0&1\cr 
&&&0&0\cr \end{matrix} \cr
& 0\cr
&&\ddots\cr
&&&0\cr
\end{matrix}\right],\quad\quad 0\le r\le\frac n2,
\]
\[
\left[
\kern-0.8cm
\begin{matrix}
{\scriptstyle 2r}\left\{
\phantom{\begin{matrix} 0\cr 0\cr \vdots\cr 0\cr 0\cr 0\cr 0\cr \end{matrix}}\right.
\begin{matrix} 
0&1\cr 
0&0\cr 
&&\ddots\cr 
&&&0&1\cr 
&&&0&0\cr 
&&&&&1&1\cr
&&&&&0&1\cr \end{matrix} \cr
& 0\cr
&&\ddots\cr
&&&0\cr
\end{matrix}\right],\quad\quad 1\le r\le\frac n2,
\]
\[
\left[
\kern-0.8cm
\begin{matrix}
{\scriptstyle 2r}\left\{
\phantom{\begin{matrix} 0\cr 0\cr \vdots\cr 0\cr 0\cr \end{matrix}}\right.
\begin{matrix} 
0&1\cr 
0&0\cr 
&&\ddots\cr 
&&&0&1\cr 
&&&0&0\cr \end{matrix} \cr
& 1\cr
&&0\cr
&&&\ddots\cr
&&&&0\cr
\end{matrix}\right],\quad\quad 0\le r\le\frac {n-1}2.
\]
We say that the matrix $A$ and its corresponding quadratic form $f=$\break 
$(x_1,\dots,x_n)A(x_1,\dots,x_n)^T$ are of type $(n,r,0,0)$ or $(n,r,1,0)$
or $(n,r,0,1)$ according to the above three types of canonical forms of $A$. In a type $(n,r,u,v)$, $n,r,v\in\Bbb N$ but $u\in\Bbb F_2$. The type of $A$ or $f$ is denoted by $\text{\rm type}(A)$ or $\text{\rm type}(f)$.

\item[(ii)]
We have
\[
\left[
\begin{matrix}
1\cr
&1\cr
\end{matrix}\right]\cong\left[
\begin{matrix}
1&0\cr
0&0\cr
\end{matrix}\right],\quad
\left[
\begin{matrix}
1&1\cr
0&1\cr
&&1\cr
\end{matrix}\right]\cong \left[
\begin{matrix}
0&1\cr
0&0\cr
&&1\cr
\end{matrix}\right],\quad
\left[
\begin{matrix}
1&1\cr
0&1\cr
&&1&1\cr
&&0&1\cr
\end{matrix}\right]\cong
\left[
\begin{matrix}
0&1\cr
0&0\cr
&&0&1\cr
&&0&0\cr
\end{matrix}\right].
\]

\item[(iii)]
If $\text{\rm type}(A_i)=(n_i,r_i,u_i,v_i)$, $i=1,2$, then
\begin{equation}\label{3.1}
\text{\rm type}(A_1\oplus A_2)=\bigl(n_1+n_2,\, r_1+r_2,\, u_1+u_2,\,\delta(v_1+v_2)\bigr).
\end{equation}

\item[(iv)]
If $\text{\rm type}(A)=(n,r,u,v)$ and $m\ge 0$, then
\begin{equation}\label{3.2}
\text{\rm type}(A\otimes I_m)=\text{\rm type}(\underbrace{A\oplus\cdots\oplus A}_m)=
\bigl(mn,\, mr,\, mu,\,\delta(m)v\bigr).
\end{equation}

\item[(v)]
Let $f(x_1,\dots,x_n)$ be a quadratic form over $\Bbb F_2$ of type $(n,r,u,v)$. Then
\begin{equation}\label{3.3}
|Z(f)|=2^{n-1}+(1-v)(-1)^u 2^{n-1-r}.
\end{equation}
\end{enumerate}
\end{thm}

In Theorem~\ref{T3.1}, (i) is Dickson's theorem \cite[Theorem 199]{Dic58}
specialized for $\Bbb F_2$; (ii) -- (v) can be verified easily; (v) is 
the main reason that binary quadratic forms are used in many areas, see, for example,
Dillon and Dobbertin \cite[Appendix A]{Dil04}. It is important to observe that if $\text{type}(f)=(n,r,0,1)$, $|Z(f)|=2^{n-1}$ is
independent of $r$.

\begin{thm}\label{TA}
Let $N(n,r,u,v)$ be the number of quadratic forms of type $(n,r,u,v)$ in $x_1,\dots, 
x_n$ over $\Bbb F_2$. Then
\begin{align*}
N(n,r,0,1)&=2^{r(r+1)}\frac{\prod_{i=0}^{2r}(2^{n-2r+i}-1)}{\prod_{i=1}^r(2^{2i}-1)},\\
N(n,r,0,0)&=2^{r^2-1}(2^r+1)\frac{\prod_{i=1}^{2r}(2^{n-2r+i}-1)}{\prod_{i=1}^r(2^{2i}-1)},\\ 
N(n,r,1,0)&=2^{r^2-1}(2^r-1)\frac{\prod_{i=1}^{2r}(2^{n-2r+i}-1)}{\prod_{i=1}^r(2^{2i}-1)}.
\end{align*}
\end{thm}

\begin{proof}
Let $R_0(s,n)$ be the set of all polynomial functions $f(x_1,\dots,x_n)$ from $\Bbb F_2^n$ to
$\Bbb F_2$ such that $\deg f\le s$ and $f(0)=0$. 
Let $N(n,r)$ be the number of elements in $R_0(2,n)/R_0(1,n)$ which are linearly equivalent to $x_1x_2+x_3x_4+\cdots+x_{2r-1}x_{2r}$. Then
\begin{equation}\label{A-1}
N(n,r)=2^{r(r-1)}\frac{\prod_{i=1}^{2r}(2^{n-2r+i}-1)}{\prod_{i=1}^r(2^{2i}-1)}.
\end{equation}
(Cf. \cite[Ch. 15, Theorem 2]{Mac81} or \cite[Lemma 2.2]{Hou94}.) With a fixed $r$ ($0\le r\le\lfloor\frac n2\rfloor$), let $N_{(n,r)}(u,v)$ be the number of $g\in R_0(1,n)$ such that
\[
\text{type}(x_1x_2+x_3x_4+\cdots+x_{2r-1}x_{2r}+g)=(n,r,u,v).
\]
Write $g=b_1x_1+\cdots+b_nx_n$. Then
\[
\begin{split}
&\text{type}(x_1x_2+x_3x_4+\cdots+x_{2r-1}x_{2r}+g)\cr
=\,&
\begin{cases}
(n,r,0,1)&\text{if}\ (b_{2r+1},\dots,b_n)\ne 0,\cr
(n,r,0,0)&\text{if}\ (b_{2r+1},\dots,b_n)= 0\ \text{and}\ b_1b_2+b_3b_4+\cdots+b_{2r-1}b_{2r}=0,\cr
(n,r,1,0)&\text{if}\ (b_{2r+1},\dots,b_n)= 0\ \text{and}\ b_1b_2+b_3b_4+\cdots+b_{2r-1}b_{2r}=1.
\end{cases}
\end{split}
\]
Therefore,
\begin{equation}\label{A-2}
\begin{cases}
N_{(n,r)}(0,1)=2^{2r}(2^{n-2r}-1),\cr
N_{(n,r)}(0,0)=2^{2r-1}+2^{r-1},\cr
N_{(n,r)}(1,0)=2^{2r-1}-2^{r-1}.
\end{cases}
\end{equation}
Since $N(n,r,u,v)=N(n,r)\cdot N_{(n,r)}(u,v)$, the conclusion of the theorem follows immediately form \eqref{A-1} and \eqref{A-2}.
\end{proof}

We can define addition and scalar multiplication for types. Let $A_1$ and $A_2$ be square matrices over $\Bbb F_2$ and let
$m\in \Bbb N$. We define
\[
\text{type}(A_1)\boxplus\text{type}(A_2)=\text{type}(A_1\oplus A_2)
\]
and
\[
m*\text{type}(A_1)=\text{type}(\underbrace{A_1\oplus\cdots\oplus A_1}_m)=\text{type}(A_1\otimes I_m)=
\underbrace{\text{type}(A_1)\boxplus\cdots\boxplus \text{type}(A_1)}_m.
\]
Then \eqref{3.1} and \eqref{3.2} become
\begin{equation}\label{3.4}
(n_1,r_1,u_1,v_1)\boxplus(n_2,r_2,u_2,v_2)
=\bigl(n_1+n_2,\, r_1+r_2,\, u_1+u_2,\, \delta(v_1+v_2)\bigr)
\end{equation}
and
\begin{equation}\label{3.5}
m*(n,r,u,v)=\bigl(mn,\,mr,\,mu,\,\delta(m)v\bigr).
\end{equation}

In the subsequent sections, we will need to determine the number of common zeros of several quadratic forms. The following lemma is useful for this purpose.

\begin{lem}\label{L3.3}
Let $f_1,\dots, f_r$ be functions from $\Bbb F_2^n$ to $\Bbb F_2$. Then
\[
|Z(f_1)\cap\cdots\cap Z(f_r)|=-2^n+\frac 1{2^{r-1}}\sum_{(a_1,\dots,a_r)\in\Bbb F_2^r}|Z(a_1f_1+\cdots+a_rf_r)|.
\]
\end{lem}

\begin{proof}
We have
\[
\begin{split}
&\sum_{(a_1,\dots,a_r)\in\Bbb F_2^r}|Z(a_1f_1+\cdots+a_rf_r)|
\vrule height0pt width0pt depth20pt \cr
=\,&\sum_{x\in\Bbb F_2^n}\sum_{\substack{(a_1,\dots,a_r)\in\Bbb F_2^r\cr
a_1f_1(x)+\cdots+a_rf_r(x)=0}} 1
\vrule height0pt width0pt depth28pt \cr
=\,& |Z(f_1)\cap\cdots\cap Z(f_r)|\cdot 2^r+\bigl(2^n-|Z(f_1)\cap\cdots\cap Z(f_r)|\bigr)2^{r-1}
\vrule height0pt width0pt depth8pt \cr
=\,& |Z(f_1)\cap\cdots\cap Z(f_r)|\cdot 2^{r-1}+2^{n+r-1}.
\end{split}
\]
It follows that
\[
|Z(f_1)\cap\cdots\cap Z(f_r)|= -2^n+\frac 1{2^{r-1}}
\sum_{(a_1,\dots,a_r)\in\Bbb F_2^r}|Z(a_1f_1+\cdots+a_rf_r)|.
\]
\end{proof}


\section{Outline of the Method and Preliminary Results}

Let 
\[
S_{k\times n}=\{X\in M_{k\times n}: XX^T=0\}.
\]
We will treat the elements in the symmetric group $\frak S_n$ as $n\times n$ permutation matrices. The group ${\rm GL}(k,\Bbb F_2)\times\frak S_n$ acts on $S_{k\times n}$ as
follows: For $A\in {\rm GL}(k,\Bbb F_2)$, $P\in\frak S_n$ and $X\in S_{k\times n}$,
\[
X^{(A,P)}=A^{-1}XP.
\]
Each $X\in S_{k\times n}$ generates a code (the row space of $X$) in $\Bbb F_2^n$ which is self-orthogonal and of dimension $\le k$. Two matrices $X_1,X_2\in S_{k\times n}$ generate equivalent codes if and only if $X_1$ and $X_2$ are in the same ${\rm GL}(k,\Bbb F_2)\times\frak S_n$-orbit. Let $\Psi_{\le k,n}$ be the number of ${\rm GL}(k,\Bbb F_2)\times\frak S_n$-orbits in 
$S_{k\times n}$. Then $\Psi_{\le k,n}$ is the number of inequivalent self-orthogonal codes in $\Bbb F_2^n$ of dimension $\le k$. Clearly,
\[
\Psi_{k,n}=\Psi_{\le k,n}-\Psi_{\le k-1, n}.
\]
Therefore, to compute $\Psi_{k,n}$, it suffices to compute $\Psi_{\le k,n}$. We will concentrate
on $\Psi_{\le k,n}$ in the paper.

We will need two notions of equivalence between matrices. For $X_1,X_2\in M_{k\times n}$, we write $X_1\approx X_2$ if there exits $P\in\frak S_n$ such that $X_1=X_2P$; we write
$X_1\sim X_2$ if there exist $A\in{\rm GL}(k,\Bbb F_2)$ and $P\in\frak S_n$ such that
$X_1=A^{-1}X_2P$.

We first take a moment to determine $\Psi_{\le 0,n}$, $\Psi_{\le 1,n}$ and $\Psi_{\le 2,n}$.
Obviously,
\[
\begin{split}
\Psi_{\le 0,n}\,&=1,\cr
\Psi_{\le 1,n}\,&=\left\lfloor\frac n2\right\rfloor+1.
\end{split}
\]

\begin{prop}\label{P2.1}
We have
\begin{equation}\label{2.1}
\begin{split}
\Psi_{\le 2,n}=\,&\frac13\Bigl(\left\lfloor\frac n6\right\rfloor+1\Bigr)+\frac 12 \Bigl(\left\lfloor\frac n2\right\rfloor-\left\lfloor\frac n4\right\rfloor+1\Bigr)\Bigl(\left\lfloor\frac n4\right\rfloor+1\Bigr)\cr
&+\frac 1{36}\Bigl(\left\lfloor\frac n2\right\rfloor+3\Bigr)\Bigl(\left\lfloor\frac n2\right\rfloor+2\Bigr)
\Bigl(\left\lfloor\frac n2\right\rfloor+1\Bigr).
\end{split}
\end{equation}
\end{prop}

\begin{proof}
Every matrix in $S_{2\times n}$ is $\sim$ equivalent to a matrix of the form 
\[
X_{a,b,c}:=\left[
\phantom{\begin{matrix} 1\cr 1 \end{matrix}}
\right.\kern-2mm
\overbrace{
\begin{matrix}
1\ \cdots\ 1\cr
1\ \cdots\ 1\cr
\end{matrix}}^a\ 
\overbrace{
\begin{matrix}
1\ \cdots\ 1\cr
0\ \cdots\ 0\cr
\end{matrix}}^b\
\overbrace{
\begin{matrix}
0\ \cdots\ 0\cr
1\ \cdots\ 1\cr
\end{matrix}}^c\ 
\begin{matrix}
0\ \cdots\ 0\cr
0\ \cdots\ 0\cr
\end{matrix}
\left.\kern-2mm
\phantom{\begin{matrix} 1\cr 1 \end{matrix}}
\right],
\]
where $a,b,c$ are all even. Moreover, $X_{a,b,c}\sim X_{a',b',c'}$ if and only if $a',b',c'$ is
a permutation of $a,b,c$. Therefore,
\[
\Psi_{\le 2,n}=\bigl|\bigl\{ (a_1,a_2,a_3)\in\Bbb N^3:0\le a_1\le a_2\le a_3,\ a_1+a_2+a_3\le
\frac n2\bigr\}\bigr|.
\]
Let
\[
\mathcal A=\bigl\{(a_1,a_2,a_3)\in\Bbb N^3:a_1+a_2+a_3\le\frac n2 \bigr\},
\]
\[
\mathcal A_{123}=\bigl\{(a,a,a)\in\mathcal A\bigr\}
\]
and for $i,j\in\{1,2,3\}$, let
\[
\mathcal A_{ij}=\{(a_1,a_2,a_3)\in\mathcal A:a_i=a_j\}.
\]
Then
\[
\begin{split}
\Psi_{\le 2,n}\,&=|\mathcal A_{123}|+\frac 13|(\mathcal A_{12}\cup
\mathcal A_{13}\cup\mathcal A_{23})\setminus\mathcal A_{123}|+
\frac 16|\mathcal A\setminus(\mathcal A_{12}\cup
\mathcal A_{13}\cup\mathcal A_{23})|\cr
&=\frac 23 |\mathcal A_{123}|+\frac 16 |\mathcal A_{12}\cup
\mathcal A_{13}\cup\mathcal A_{23}|+\frac 16|\mathcal A|.
\end{split}
\]
By the inclusion-exclusion formula,
\[
|\mathcal A_{12}\cup
\mathcal A_{13}\cup\mathcal A_{23}|=3|\mathcal A_{12}|-2|\mathcal A_{123}|.
\]
Hence 
\begin{equation}\label{2.2}
\Psi_{\le 2,n}=\frac 13 |\mathcal A_{123}|+\frac 12|\mathcal A_{12}|+\frac 16 |\mathcal A|.
\end{equation} 
Clearly, 
\begin{equation}\label{2.3}
|\mathcal A_{123}|=\left\lfloor\frac n6\right\rfloor+1,
\end{equation}
\begin{equation}\label{2.4}
|\mathcal A|=\bigl|\bigl\{(a_1,a_2,a_3,a_4)\in\Bbb N^4: a_1+a_2+a_3+a_4=\left\lfloor
\frac n2\right\rfloor\bigr\}\bigr|=\binom{\lfloor \frac n2\rfloor+3}3
\end{equation}
and
\begin{equation}\label{2.5}
\begin{split}
|\mathcal A_{12}|\,&=\sum_{\substack{a_1,a_3\in\Bbb N\cr
2a_1+a_3\le\frac n2}}1  \vrule height0pt width0pt depth25pt \cr
&=\sum_{0\le a_1\le\lfloor\frac n4\rfloor}\ \sum_{0\le a_3\le\lfloor \frac n2\rfloor-2a_1}1
\vrule height0pt width0pt depth20pt \cr
&=\sum_{0\le a_1\le\lfloor\frac n4\rfloor}\Bigl(\left\lfloor\frac n2\right\rfloor-2a_1+1\Bigr)
\vrule height0pt width0pt depth20pt \cr
&=\Bigl(\left\lfloor\frac n2\right\rfloor+1\Bigr)\Bigl(\left\lfloor\frac n4\right\rfloor+1\Bigr)
-\left\lfloor\frac n4\right\rfloor\Bigl(\left\lfloor\frac n4\right\rfloor+1\Bigr)
\vrule height0pt width0pt depth15pt \cr
&=\Bigl(\left\lfloor\frac n2\right\rfloor-\left\lfloor\frac n4\right\rfloor+1\Bigr)
\Bigl(\left\lfloor\frac n4\right\rfloor+1\Bigr).
\end{split}
\end{equation}
Equation \eqref{2.1} follows from \eqref{2.2} -- \eqref{2.5}.
\end{proof}

Unfortunately, the combinatorial method in Proposition~\ref{P2.1} does not seem to have a 
generalization for $\Psi_{\le k,n}$ with $k\ge 3$. To compute $\Psi_{\le k,n}$ with $k\ge 3$, we start afresh with a more algebraic approach.

By the Burnside lemma,
\begin{equation}\label{2.6}
\Psi_{\le k,n}=\frac 1{|{\rm GL}(k,\Bbb F_2)\times\frak S_n|}\sum_{\substack{A\in
{\rm GL}(k,\Bbb F_2)\cr P\in\frak S_n}}|{\rm Fix}(A,P)|,
\end{equation}
where
\[
\begin{split}
{\rm Fix}(A,P)\,&=\{X\in S_{k\times n}:X^{(A,P)}=X\}\cr
&=\{X\in M_{k\times n}:AX=XP,\ XX^T=0\}.
\end{split}
\]

If $G$ is a group, $\mathcal C(G)$ denotes a set of representatives of the conjugacy classes
of $G$. For $g\in G$, ${\rm cent}_G(g)$ denotes the centralizer of $g$ in $G$. A partition of an integer $n>0$ is a sequence of nonnegative integers $\lambda=(\lambda_1,\lambda_2,\dots)$ such
that $\sum_{i\ge 1}i\lambda_i=n$. We write $\lambda\vdash n$ to mean that $\lambda$ is a partition of $n$. 
For $\lambda=(\lambda_1,\lambda_2,\dots)\vdash n$ and $0\le a<b$, we define
\[
\lambda_{a,b}=\sum_{i\equiv a\, (b)}\lambda_i.
\]
This expression will appear repreatedly in formulas later on. 
For each $\lambda\vdash n$, let $P_\lambda\in\frak S_n$ be the ``canonical''
permutation of cycle type $\lambda$. For example, if $\lambda=(0,2,1)\vdash 7$, then $P_\lambda=(1,2)(3,4)(5,6,7)$ as a permutation and
\[
P_\lambda=\left[
\begin{matrix}
0&1\cr
1&0\cr
\end{matrix}\right]\oplus
\left[
\begin{matrix}
0&1\cr
1&0\cr
\end{matrix}\right]\oplus
\left[
\begin{matrix}
0&1&0\cr
0&0&1\cr
1&0&0\cr
\end{matrix}\right]
\]
as a permutation matrix. We can choose $\mathcal C(\frak S_n)=\{P_\lambda:\lambda\vdash n\}$.
If $\lambda=(\lambda_1,\lambda_2,\dots)\vdash n$, it is well known that
\[
|{\rm cent}_{\frak S_n}(P_\lambda)|=\lambda_1!\lambda_2!\cdots 1^{\lambda_1}2^{\lambda_2}\cdots.
\]
Therefore, we can write \eqref{2.6} as 
\begin{equation}\label{2.7}
\begin{split}
\Psi_{\le k,n}\,&=\sum_{\substack{A\in\mathcal C({\rm GL}(k,\Bbb F_2))\cr
P\in\mathcal C(\frak S_n)}}\frac 1{|{\rm cent}_{{\rm GL}(k,\Bbb F_2)}(A)|\,
|{\rm cent}_{\frak S_n}(P)|} |{\rm Fix}(A,P)|\cr
&=\sum_{A\in\mathcal C({\rm GL}(k,\Bbb F_2))}
\frac 1{|{\rm cent}_{{\rm GL}(k,\Bbb F_2)}(A)|}\sum_{\lambda=(\lambda_1,\lambda_2,\dots)\vdash n}
\frac{|{\rm Fix}(A,P_\lambda)|}{\lambda_1!\lambda_2!\cdots 1^{\lambda_1}2^{\lambda_2}\cdots}.
\end{split}
\end{equation}
In \eqref{2.7}, the elements in $\mathcal C({\rm GL}(k,\Bbb F_2))$ are the canonical forms of $k\times k$ invertible matrices and they can be enumerated in terms of their elementary divisors. As for 
$|{\rm cent}_{{\rm GL}(k,\Bbb F_2)}(A)|$, it suffices to assume that all the elementary divisors of $A$ are powers of a single irreducible polynomial
in $\Bbb F_2[x]$ since if $A=A_1\oplus A_2$, $A_i\in{\rm GL}(k_i,\Bbb F_2)$, where every elementary divisor of $A_1$ is prime to every elementary
divisor of $A_2$, then
\[
|{\rm cent}_{{\rm GL}(k,\Bbb F_2)}(A)|=
|{\rm cent}_{{\rm GL}(k_1,\Bbb F_2)}(A_1)|\,
|{\rm cent}_{{\rm GL}(k_2,\Bbb F_2)}(A_2)|.
\]

\begin{thm}\label{T2.2} 
{\rm (\cite[Theorem 3.6]{Hou96})} 
Assume that $A$ is a $k\times k$ matrix over $\Bbb F_q$ with elementary
divisors $\underbrace{f^1,\dots, f^1}_{\mu_1}, 
\underbrace{f^2,\dots, f^2}_{\mu_2},\dots$, where $f\in\Bbb F_q[x]$ is irreducible of degree $d$. Then
\[
|{\rm cent}_{{\rm GL}(k,\Bbb F_q)}(A)|
=\prod_{i\ge 1}q^{d\mu_i(1\mu_1+2\mu_2+\cdots+i\mu_i+i\mu_{i+1}+\cdots)}
\prod_{j=1}^{\mu_i}(1-q^{-dj}).
\]
\end{thm}

Now, the only component in \eqref{2.7} that needs to be determined is
$|{\rm Fix}(A,P_\lambda)|$.

\begin{thm}\label{T2.3}
Let $\lambda=(\lambda_1,\lambda_2,\dots)\vdash n$ and
let $A\in{\rm GL}(k,\Bbb F_2)$ with multiplicative order $o(A)=t$.
For each $d\mid t$, let $s_d=k-{\rm rank}(A^d-I)$, let $B_d\in M_{k\times s_d}$ such that its columns form a basis of 
\[
\{x\in\Bbb F_2^k:(A^d-I)x=0\},
\]
and let
\begin{equation}\label{2.8}
\alpha_d=\sum_{\substack{i\ge 1,\, \nu(i)\le \nu(t) \cr {\rm gcd}(i,t)=d}}
\lambda_i.
\end{equation}           
Then 
\begin{equation}\label{2.9}
|{\rm Fix}(A,P_\lambda)|=2^{\sum_{\nu(i)>\nu(t)} s_{{\rm gcd}(i,t)}\lambda_i}
\cdot\frak n(A),
\end{equation}
where $\frak n(A)$ is the number of sequences of matrices $(Y_d)_{d\mid t}$ with $Y_d\in M_{s_d\times \alpha_d}$ and
\begin{equation}\label{2.10}
\sum_{d\mid t}\sum_{j=0}^{d-1}
A^j B_dY_dY_d^T B_d^T(A^j)^T=0.
\end{equation}
\end{thm}

\begin{proof}
First note that $|{\rm Fix}(A,P_\lambda)|=|{\rm Fix}(A,P_\lambda^{-1})|$
since $P_\lambda$ and $P_\lambda^{-1}$ are conjugates. We will compute 
$|{\rm Fix}(A,P_\lambda^{-1})|$ since the notation is more convenient for
$|{\rm Fix}(A,P_\lambda^{-1})|$.

Let $X\in M_{k\times n}$. Then $AX=XP_\lambda^{-1}$ if and only if
\begin{equation}\label{2.11}
X=\bigl[\dots;\ 
\underbrace{x_1^{(i)},Ax_1^{(i)},\dots,A^{i-1}x_1^{(i)}}_i;\
\dots;
\underbrace{x_{\lambda_i}^{(i)},Ax_{\lambda_i}^{(i)},\dots,A^{i-1}x_{\lambda_i}^{(i)}}_i;\ \dots\bigr],
\end{equation}
\vskip -7mm
\[
\kern1.5cm\underbrace{\kern4.2cm}_{\lambda_i}
\]
where $i\ge 1$ and $x_j^{(i)}\in\Bbb F_2^k$ with $(A^i-I)x_j^{(i)}=0$, $1\le j\le\lambda_i$. 
Put $d_i={\rm gcd}(i,t)$. Note that $x_1^{(i)},\dots, x_{\lambda_i}^{(i)}$ are in the column space of $B_{d_i}$, hence
$[x_1^{(i)},\dots, x_{\lambda_i}^{(i)}]=B_{d_i}X_i$ for some $X_i\in M_{s_{d_i}\times \lambda_i}$. Therefore 
we can write \eqref{2.11} as
\begin{equation}\label{2.12}
X\approx\bigl[\dots;\ B_{d_i}X_i, AB_{d_i}X_i, \dots, A^{i-1}B_{d_i}X_i;\
\dots\bigr],
\end{equation}
where $X_i\in M_{s_{d_i}\times\lambda_i}$ is arbitrary but the equivalence
$\approx$ is given by a fixed column permutation. Since $(A^{d_i}-I)B_{d_i}=0$, i.e., $A^{d_i}B_{d_i}=B_{d_i}$, we have
\begin{equation}\label{2.13}
\bigl[\ B_{d_i}X_i, AB_{d_i}X_i, \dots, A^{i-1}B_{d_i}X_i\bigr]\approx 
\bigl[\ B_{d_i}X_i, AB_{d_i}X_i, \dots, A^{d_i-1}B_{d_i}X_i\bigr]\otimes
{\bf 1}(\frac i{d_i}).
\end{equation}
From \eqref{2.12} and \eqref{2.13}, we have
\[
\begin{split}
XX^T\,&=\sum_{\substack{i\ge 1\cr i/d_i\ \text{odd}}}\sum_{j=0}^{d_i-1}
A^jB_{d_i}X_iX_i^TB_{d_i}^T(A^j)^T\cr
&=\sum_{d\mid t}\sum_{j=0}^{d-1}
A^jB_dY_dY_d^T B_d^T(A^j)^T,
\end{split}
\]
where $Y_d$ is the concatenation of all $X_i$ with $\nu(i)\le\nu(t)$ and
${\rm gcd}(i,t)=d$. Thus $XX^T=0$ if and only if \eqref{2.10} holds. The equation $XX^T=0$ imposes no restriction on $X_i$ with $\nu(i)>\nu(t)$
and the number of $X_i$ with $\nu(i)>\nu(t)$ is 
\[
2^{\sum_{\nu(i)>\nu(t)}s_{d_i}\lambda_i}.
\]
Therefore, we have
\[
\begin{split}
|{\rm Fix}(A,P_\lambda^{-1})|\,&= 2^{\sum_{\nu(i)>\nu(t)}s_{d_i}\lambda_i}
\cdot\bigl(\text{the number of $(Y_d)_{d\mid t}$
satisfying \eqref{2.10}}\bigr)\cr
&=2^{\sum_{\nu(i)>\nu(t)}s_{d_i}\lambda_i}\cdot\frak n(A).
\end{split}
\]
\end{proof}

In Theorem~\ref{T2.3} (equation~\eqref{2.9}), the remaining question is how to compute the number $\frak n(A)$. We now set to answer this question.

Put $Y_dY_d^T=[z_{ij}^{(d)}]_{s_d\times s_d}$ where $z_{ij}^{(d)}=z_{ji}^{(d)}$. Then each entry of the left side of \eqref{2.10} is a linear
function of $z_{ij}^{(d)}$, $0\le i\le j\le s_d$, $d\mid t$. Thus \eqref{2.10} can be written as a system of $L$ equations
\begin{equation}\label{2.12.1}
\sum_{d\mid t}\sum_{0\le i\le j\le s_d}c_{l;ij}^{(d)} z_{ij}^{(d)}=0,\qquad 1\le l\le L,
\end{equation}
where $c_{l;ij}^{(d)}\in\Bbb F_2$ are constants. We define $c_{l;ij}^{(d)}=0$ for $i>j$ and let $C_l^{(d)}=[c_{l;ij}^{(d)}]_{s_d\times s_d}$.
Although the left side of \eqref{2.10} is a symmetric matrix of size $k\times k$, the number $L$ of equations in \eqref{2.12.1} is usually much smaller
than $\frac 12 k(k+1)$, as we will see in actual computations. 

Write $Y_d=[y_{ij}^{(d)}]_{s_d\times \alpha_d}$ and let $\widetilde Y_d$ be the concatenation of all the rows of $Y_d$, i.e.,
\[
\widetilde Y_d=\bigl[ y_{11}^{(d)},\dots,y_{1,\alpha_d}^{(d)},\ \dots,\ y_{s_d,1}^{(d)},\dots,y_{s_d,\alpha_d}^{(d)}\bigr]_{1\times(s_d\alpha_d)}.
\]
Since 
\[
\bigl[ y_{i1}^{(d)},\dots,y_{i,\alpha_d}^{(d)}\bigr]\bigl[ y_{j1}^{(d)},\dots,y_{j,\alpha_d}^{(d)}\bigr]^T=z_{ij}^{(d)},
\]
we have
\[
\sum_{0\le i\le j\le s_d}c_{l;ij}^{(d)}z_{ij}^{(d)}=\widetilde Y_d(C_l^{(d)}\otimes I_{\alpha_d})\widetilde Y_d^T. 
\]
Let $Y_{1\times\theta}$ be the concatenation of all $\widetilde Y_d$, $d\mid t$, where 
\begin{equation}\label{I1.1}
\theta=\sum_{d\mid t}s_d\alpha_d.
\end{equation}
Then \eqref{2.12.1} becomes 
\begin{equation}\label{2.12.2}
Y\Bigl[\bigoplus_{d\mid t}\bigl(C_l^{(d)}\otimes I_{\alpha_d}\bigr)\Bigr]Y^T=0,\qquad 1\le l\le L.
\end{equation}
The left side of \eqref{2.12.2} is a quadratic form in $Y$, which is denoted by $f_l(Y)$. Therefore,
\begin{equation}\label{2.12.3}
\begin{split}
\frak n(A)\,&=|Z(f_1)\cap\cdots\cap Z(f_L)|\cr
&=-2^\theta+2^{-L+1}\sum_{(a_1,\dots,a_L)\in\Bbb F_2^L}|Z(a_1f_1+\cdots+a_Lf_L)|\qquad\text{(by Lemma~\ref{L3.3})}.
\end{split}
\end{equation}
In \eqref{2.12.3}, we have
\[
a_1f_1+\cdots+a_Lf_L=Y\biggl[\bigoplus_{d\mid t}\Bigl[\Bigl(\sum_{l=1}^L a_lC_l^{(d)}\Bigr)\otimes I_{\alpha_d}\Bigr]\biggr]Y^T.
\] 
Put
\[
\frak C_{(a_1,\dots,a_L)}=\bigoplus_{d\mid t}\Bigl[\Bigl(\sum_{l=1}^L a_lC_l^{(d)}\Bigr)\otimes I_{\alpha_d}\Bigr],\qquad
(a_1,\dots,a_L)\in\Bbb F_2^L.
\]
By Theorem~\ref{T3.1}, $|Z(a_1f_1+\cdots+a_Lf_L)|=|Z(Y\frak C_{(a_1,\dots,a_L)}Y^T)|$ is determined by
\[
\text{type}(\frak C_{(a_1,\dots,a_L)})=\underset{d\mid t}{\vcenter{\hbox{\huge $\boxplus$}}}
\Bigl[\alpha_d*\text{type}\Bigl(\sum_{l=1}^La_lC_l^{(d)}\Bigr)\Bigr].
\]

To sum up, we have the following algorithm for computing $\frak n(A)$.

\begin{algo}[An algorithm for computing $\frak n(A)$]\label{A1}\
\begin{list}%
{\thestep.}
{\usecounter{step}   
\setlength {\leftmargin 1.5cm}
\setlength{\labelwidth 1.3cm}
}

\item
Let $Y_dY_d^T=[z_{ij}^{(d)}]$ and write \eqref{2.10} entry wise, in the form of 
\[
\sum_{d\mid t}\sum_{0\le i\le j\le s_d}c_{l;ij}^{(d)}z_{ij}^{(d)}=0,\qquad 1\le l\le L.
\]
Record the matrices $C_l^{(d)}=[c_{l;ij}^{(d)}]$, $1\le l\le L$, $d\mid t$.

\item 
For each $(a_1,\dots,a_L)\in\Bbb F_2^L$ and $d\mid t$, determine $\text{\rm type}\bigl(\sum_{l=1}^L a_lC_l^{(d)}\bigr)$.

\item
For each $(a_1,\dots,a_L)\in\Bbb F_2^L$, put $\frak C_{(a_1,\dots,a_L)}=\bigoplus_{d\mid t}\bigl[\bigl(\sum_{l=1}^L a_lC_l^{(d)}\bigr)\otimes I_{\alpha_d}\bigr]$ and compute 
\[
\text{\rm type}(\frak C_{(a_1,\dots,a_L)})=\underset{d\mid t}{\vcenter{\hbox{\huge $\boxplus$}}}
\Bigl[\alpha_d*\text{\rm type}\Bigl(\sum_{l=1}^La_lC_l^{(d)}\Bigr)\Bigr].
\]
Find $|Z(Y\frak C_{(a_1,\dots,a_L)}Y^T)|$ by \eqref{3.3}.

\item
By \eqref{2.12.3}
\begin{equation}\label{n(A)}
\frak n(A)=-2^\theta+2^{-L+1}\sum_{(a_1,\dots,a_L)\in\Bbb F_2^L}|Z(Y\frak C_{(a_1,\dots,a_L)}Y^T)|.
\end{equation}
\end{list}
\end{algo}

When $k$ is small, the above algorithm is effective, as we will see in the following examples and in the computations of the subsequent sections.

\begin{exmp}\label{E1}
\rm
Let $A=\genfrac{[}{]}{0pt}{1}{0\ 1}{1\ 0} \in{\rm GL}(2,\Bbb F_2)$. We determine $\frak n(A)$ by Algorithm~\ref{A1}.
We have $o(A)=2$.  Since
\[
A-I=
\left[
\begin{matrix}
1&1\cr
1&1\cr
\end{matrix}\right]\quad\text{and}\quad  
A^2-I=0,
\]
we have $s_1=1$, $s_2=2$. We can choose
\[
B_1=\left[
\begin{matrix}
1\cr
1\cr
\end{matrix}\right],\quad B_2=I_2.
\]
By \eqref{2.8} and \eqref{I1.1},
\begin{equation}\label{2.12.4}
\alpha_1=\sum_{\nu(i)=0}\lambda_i=\lambda_{1,2},\quad
\alpha_2=\sum_{\nu(i)=1}\lambda_i=\lambda_{2,4},\quad \theta=\alpha_1+2\alpha_2.
\end{equation}
Now we are in Step 1 of Algorithm~\ref{A1}. Put $Y_dY_d^T=[z_{ij}^{(d)}]$. Then \eqref{2.10} becomes
\[
\left[
\begin{matrix}
z_{11}^{(1)}&z_{11}^{(1)}\cr
z_{11}^{(1)}&z_{11}^{(1)}\cr
\end{matrix}\right]+
\left[
\begin{matrix}
z_{11}^{(2)}+z_{22}^{(2)}&0\cr
0&z_{11}^{(2)}+z_{22}^{(2)}\cr
\end{matrix}\right]=0,
\]
i.e.,
\[
\begin{cases}
z_{11}^{(1)}=0,\cr
z_{11}^{(2)}+z_{22}^{(2)}=0.
\end{cases}
\]
Thus, $L=2$ and 
\begin{align*}
C_1^{(1)}&=[1], &C_1^{(2)}&=0,\\
C_2^{(1)}&=[0], &C_2^{(2)}&=I_2.
\end{align*}

In Step 2 of Algorithm~\ref{A1}, we find that
\[
\text{type}(a_1C_1^{(1)}+a_2C_2^{(1)})=
\begin{cases}
(1,0,0,0)&\text{if}\ a_1=0,\cr
(1,0,0,1)&\text{if}\ a_1=1,
\end{cases}
\]
\[
\text{type}(a_1C_1^{(2)}+a_2C_2^{(2)})=
\begin{cases}
(2,0,0,0)&\text{if}\ a_2=0,\cr
(2,0,0,1)&\text{if}\ a_2=1.
\end{cases}
\]

In Step 3, we have
\[
\begin{split}
\text{type}(\frak C_{(a_1,a_2)})\,&=\alpha_1*\text{type}(a_1C_1^{(1)}+a_2C_2^{(1)})\boxplus
\alpha_2*\text{type}(a_1C_1^{(2)}+a_2C_2^{(2)})\cr
&=
\begin{cases}
(\theta,0,0,0)&\text{if}\ (a_1,a_2)=(0,0),\cr
(\theta,0,0,\delta(\alpha_1))&\text{if}\ (a_1,a_2)=(1,0),\cr
(\theta,0,0,\delta(\alpha_2))&\text{if}\ (a_1,a_2)=(0,1),\cr
(\theta,0,0,\delta(\alpha_1+\alpha_2))&\text{if}\ (a_1,a_2)=(1,1)
\end{cases}
\end{split}
\]
and 
\[
|Z(Y\frak C_{(a_1,a_2)}Y^T)|=
\begin{cases}
2^\theta&\text{if}\ (a_1,a_2)=(0,0),\cr
(2-\delta(\alpha_1))2^{\theta-1}&\text{if}\ (a_1,a_2)=(1,0),\cr
(2-\delta(\alpha_2))2^{\theta-1}&\text{if}\ (a_1,a_2)=(0,1),\cr
(2-\delta(\alpha_1+\alpha_2))2^{\theta-1}&\text{if}\ (a_1,a_2)=(1,1).
\end{cases}
\]

Finally, in Step 4, we arrive at
\[
\begin{split}
\frak n(A)\,&=-2^\theta + 2^{-1}\cdot 2^{\theta-1}\bigl[ 8-\delta(\alpha_1)-\delta(\alpha_2)-\delta(\alpha_1+\alpha_2)\bigr]\cr
&=2^{\alpha_1+2\alpha_2-2}\bigl[ 4-\delta(\alpha_1)-\delta(\alpha_2)-\delta(\alpha_1+\alpha_2)\bigr].
\end{split}
\]
By \eqref{2.9} and \eqref{2.12.4}, 
\begin{equation}\label{I4}
\begin{split}
\text{Fix}\bigl(\genfrac{[}{]}{0pt}{1}{0\ 1}{1\ 0},\, P_\lambda\bigr)&=2^{2\sum_{\nu(i)>1}\lambda_i}\,\frak n(A)\cr
&=2^{\lambda_{1,2}+2\lambda_{0,2}-2}\bigl[ 4-\delta(\lambda_{1,2})-\delta(\lambda_{2,4})-\delta(\lambda_{1,2}+\lambda_{2,4})\bigr].
\end{split}
\end{equation}
\end{exmp}

\begin{exmp}\label{E2}
\rm
Let $A=\genfrac{[}{]}{0pt}{1}{0\ 1}{1\ 1} \in{\rm GL}(2,\Bbb F_2)$.
In this cases, $o(A)=3$ and
\[
A^d-I
\begin{cases}
\text{is invertible}&\text{if}\ d=1,\cr
=0&\text{if}\ d=3.
\end{cases}
\]
So, $s_1=0$, $s_3=2$, and we can choose
\[
B_1=\emptyset,\quad B_3=I_2.
\]
(Note: $B_1$ is a $2\times 0$ matrix. We denote an $n\times 0$ matrix by $\emptyset$.) By \eqref{2.8},
\begin{equation}\label{2.12.5}
\alpha_3=\sum_{i\equiv 3\,(6)}\lambda_i =\lambda_{3,6}.
\end{equation}
With $Y_dY_d^T=[z_{ij}^{(d)}]$, \eqref{2.10} becomes
\[
\left[
\begin{matrix}
0&z_{11}^{(3)}+z_{12}^{(3)}+z_{22}^{(3)}\cr
z_{11}^{(3)}+z_{12}^{(3)}+z_{22}^{(3)}&0\cr
\end{matrix}\right]=0,
\]
i.e.,
\[
z_{11}^{(3)}+z_{12}^{(3)}+z_{22}^{(3)}=0.
\]
Thus, $L=1$ and
\[
C_1^{(3)}=\left[
\begin{matrix}
1&1\cr
0&1
\end{matrix}\right].
\]
We have
\[
\text{type}(a_1C_1^{(3)})=
\begin{cases}
(2,0,0,0)&\text{if}\ a_1=0,\cr
(2,1,1,0)&\text{if}\ a_1=1,
\end{cases}
\]
\[
\text{type}(\frak C_{a_1})=\alpha_3*\text{type}(a_1C_1^{(3)})=
\begin{cases}
(2\alpha_3,0,0,0)&\text{if}\ a_1=0,\cr
(2\alpha_3,\alpha_3,\alpha_3,0)&\text{if}\ a_1=1
\end{cases}
\]
and
\[
|Z(Y\frak C_{a_1}Y^T)|=
\begin{cases}
2^{2\alpha_3}&\text{if}\ a_1=0,\cr
2^{2\alpha_3-1}+(-1)^{\alpha_3}2^{\alpha_3-1}&\text{if}\ a_1=1.
\end{cases}
\]
Therefore,
\[
\frak n(A)=-2^{2\alpha_3}+2^{2\alpha_3}+2^{2\alpha_3-1}+(-1)^{\alpha_3}2^{\alpha_3-1}=2^{2\alpha_3-1}+(-1)^{\alpha_3}2^{\alpha_3-1}.
\]
By \eqref{2.9} and \eqref{2.12.5},
\begin{equation}\label{I5}
\begin{split}
\text{Fix}\bigl(\genfrac{[}{]}{0pt}{1}{0\ 1}{1\ 1},\, P_\lambda\bigr)\,& =2^{2\sum_{i\equiv 0\,(6)}\lambda_i}\,\frak n(A)\cr
&=2^{2\lambda_{0,6}-1}\bigl[2^{2\lambda_{3,6}}+(-1)^{\lambda_{3,6}}2^{\lambda_{3,6}}\bigr].
\end{split}
\end{equation}
\end{exmp}

For the rest of this section, we collect a few lemmas which will simplify the computation of $|\text{Fix}(A,P_\lambda)|$ in many cases.

\begin{lem}\label{L2.4}
We have
\begin{equation}\label{2.14}
|\text{\rm Fix}(I_k,P_\lambda)|=2^{k\lambda_{0,2}} \Bigl|S_{k\times\lambda_{1,2}}\Bigr|.
\end{equation}
\end{lem}

\begin{proof}
Clearly, $o(I_k)=1$, 
$s_1=k$ and we can choose $B_1=I_k$. By \eqref{2.8},
$\alpha_1=\sum_{\nu(i)=0}\lambda_i=\lambda_{1,2}$.
Equation \eqref{2.10} becomes
\[
Y_1Y_1^T=0,
\]
where $Y_1\in M_{k\times\alpha_1}$. The number of such $Y_1$, i.e. $\frak n(I_k)$, 
is $|S_{k\times\alpha_1}|$. Hence by \eqref{2.9},
\[
|{\rm Fix}(I_k,P_\lambda)|=2^{k\sum_{\nu(i)>0}\lambda_i}
\bigl|S_{k\times\alpha_1}\bigr|.
\]
\end{proof}

In \eqref{2.14}, $\bigl|S_{k\times\lambda_{1,2}}\bigr|$
is given by the next lemma.

\begin{lem}\label{L2.5} 
We have
\begin{equation}\label{2.15}
|S_{k\times n}|=\sum_{0\le l\le\min\{k,\frac n2\}}|\Phi_{n,l}|\, 2^{\frac12
l(l-1)}\prod_{j=k-l+1}^k(2^j-1).
\end{equation}
\end{lem}

\begin{proof}
For each $0\le l\le\min\{k,\frac n2\}$, let
\[
S_{k\times n}^{(l)}=\{X\in S_{k\times n}:{\rm rank}\, X=l\}.
\]
Put
\[
\mathcal X=\{
(X,C):X\in S_{k\times n}^{(l)},\ C\in\Phi_{n,l},\ C\ \text{is the row space
of}\ X\}.
\]
Counting the number of elements $(X,C)\in\mathcal X$ in the order of $X,C$ and in the order of $C,X$, we have
\[
|S_{k\times n}^{(l)}|=|\Phi_{n,l}|\, |\{A\in M_{k\times l}:{\rm rank}\, A=l\}|.
\]
In fact, for each $X\in S_{k\times n}^{(l)}$, there is a unique $C\in 
\Phi_{n,l}$ such that $(X,C)\in\mathcal X$. On the other hand, for each $C\in\Phi_{n,l}$ with a basis $c_1,\dots, c_l$, $(X,C)\in\mathcal X$ if and only if
\[
X=A\left[
\begin{matrix}
c_1\cr
\vdots\cr
c_l\cr
\end{matrix}\right]
\]
for some $A\in  M_{k\times l}$ with ${\rm rank}\, A=l$. Thus
\[
\begin{split}
|S_{k\times n}^{(l)}|\,&=|\Phi_{n,l}|(2^k-2^0)(2^k-2^1)\cdots(2^k-2^{l-1})\cr
&=|\Phi_{n,l}|\, 2^{\frac 12 l(l-1)}\prod_{j=k-l+1}^k(2^j-1).
\end{split}
\]
It follows that
\[
\begin{split}
|S_{k\times n}|\,&=\sum_{0\le l\le\min\{k,\frac n2\}}
|S_{k\times n}^{(l)}|\cr
&=\sum_{0\le l\le\min\{k,\frac n2\}}|\Phi_{n,l}|\, 2^{\frac 12 l(l-1)}\prod_{j=k-l+1}^k(2^j-1).
\end{split}
\]
\end{proof}

We take another look of Lemmas~\ref{L2.4} and \ref{L2.5}. Let $\alpha_1=\lambda_{1,2}$. 
Then $|S_{k\times\alpha_1}|=\frak n(I_k)$ and Lemma~\ref{L2.5} gives a formula for this number
in terms of $|\Phi_{n,l}|$. However, $\frak n(I_k)$ can also be computed directly using
Algorithm~\ref{A1}, resulting in a formula for $\frak n(I_k)$ not involving $|\Phi_{n,l}|$.  
We follow the notation of Algorithm~\ref{A1}. With $Y_1Y_1^T=[z_{ij}^{(1)}]_{k\times k}$,
\eqref{2.10} becomes
\[
z_{ij}^{(1)}=0,\qquad 1\le i\le j\le k.
\]
So, $L=\frac 12 k(k+1)$ and 
\[
\begin{split}
&C_1^{(1)}=\left[
\begin{matrix}
1&0&\cdots&0\cr
0&0&\cdots&0\cr
\vdots&\vdots&\ddots&\vdots\cr
0&0&\cdots&0
\end{matrix}\right],\ 
C_2^{(1)}=\left[
\begin{matrix}
0&1&\cdots&0\cr
0&0&\cdots&0\cr
\vdots&\vdots&\ddots&\vdots\cr
0&0&\cdots&0
\end{matrix}\right],\ \dots,\
C_{k+1}^{(1)}=\left[
\begin{matrix}
0&0&\cdots&0\cr
0&1&\cdots&0\cr
\vdots&\vdots&\ddots&\vdots\cr
0&0&\cdots&0
\end{matrix}\right],\cr 
&\dots,\ 
C_{\frac12 k(k+1)}^{(1)}=\left[
\begin{matrix}
0&0&\cdots&0\cr
0&0&\cdots&0\cr
\vdots&\vdots&\ddots&\vdots\cr
0&0&\cdots&1
\end{matrix}\right].
\end{split} 
\]
As $(a_1,\dots,a_{\frac12 k(k+1)})$ runs through $\Bbb F_2^{\frac 12k(k+1)}$,
\[
\text{type}(a_1C_1^{(1)}+\cdots+a_{\frac 12 k(k+1)}C_{\frac 12 k(k+1)}^{(1)})=(k,r,u,v)\qquad N(k,r,u,v)\ \text{times},
\]
where $N(k,r,u,v)$ is given by Theorem~\ref{TA}. Thus,
\[
\text{type}\bigl(\frak C_{(a_1,\dots,a_{\frac 12 k(k+1)})}\bigr)= \bigl(k\alpha_1,\, r\alpha_1,\, u\alpha_1,\, \delta(\alpha_1)v\bigr)\qquad N(k,r,u,v)\ \text{times}.
\]
Therefore,
\begin{equation}\label{2.24.0}
\begin{split}
&|S_{k\times\alpha_1}|\cr
=&\frak n(I_k)\cr
=&-2^{k\alpha_1}+2^{-\frac12 k(k+1)+1}
\Bigl[\sum_{0\le r\le \frac k2}N(k,r,0,0)
(2^{k\alpha_1-1}+2^{k\alpha_1-1-r\alpha_1})\cr
&+\sum_{1\le r\le \frac k2}N(k,r,1,0) \bigl(2^{k\alpha_1-1}+
(-1)^{\alpha_1}2^{k\alpha_1-1-r\alpha_1} \bigr)\cr
&+\sum_{0\le r\le \frac {k-1}2}N(k,r,0,1) \bigl(2^{k\alpha_1-1}+
(1-\delta(\alpha_1))2^{k\alpha_1-1} \bigr)\Bigr]\cr
=&-2^{k\alpha_1}+2^{-\frac12 (k+2)(k-1)} \sum_{0\le r\le \frac k2}\frac{\prod_{i=1}^{2r}(2^{k-2r+i}-1)}{\prod_{i=1}^r(2^{2i}-1)} \, 2^{r^2-1} \cr
&\cdot \Bigl[2^{k\alpha_1+r}+2^{k\alpha_1-1-r\alpha_1}\bigl(2^r+1+(-1)^{\alpha_1}(2^r-1)\bigr)\cr
&+2^{r+1}(2^{k-2r}-1)\bigl(2^{k\alpha_1}-\delta(\alpha_1)2^{k\alpha_1-1}\bigr)\Bigr].
\end{split}
\end{equation}

\begin{lem}\label{L2.8}
Let $P\in\frak S_n$, $A\in{\rm GL}(k,\Bbb F_2)$ and $B\in
{\rm GL}(l,\Bbb F_2)$. Assume that the characteristic polynomials of $A$ and $B^{-1}$ are relatively prime. Then
\begin{equation}\label{2.24}
{\rm Fix}(A\oplus B, P)=\Bigl\{\left[
\begin{matrix}
X\cr
Y\cr
\end{matrix}\right]: X\in{\rm Fix}(A,P),\ Y\in{\rm Fix}(B,P)\Bigr\}.
\end{equation}
In particular,
\begin{equation}\label{2.25}
|{\rm Fix}(A\oplus B, P)|=|{\rm Fix}(A,P)|\,
|{\rm Fix}(B,P)|.
\end{equation}
\end{lem}

\begin{proof}
It suffices to prove that in \eqref{2.24}, the right hand side is contained
in the left hand side. Let $X\in{\rm Fix}(A,P)$ and  $Y\in{\rm Fix}(B,P)$. We have $XX^T=0$, $YY^T=0$, and
\begin{equation}\label{2.26}
\left[
\begin{matrix}
A\cr
&B\cr
\end{matrix}\right]\left[
\begin{matrix}
X\cr
Y\cr
\end{matrix}\right]=\left[
\begin{matrix}
X\cr
Y\cr
\end{matrix}\right]P.
\end{equation}
It follows that
\[
XY^T=(XP)(YP)^T=(AX)(BY)^T=AXY^TB^T,
\]
i.e.,
\[
A(XY^T)=(XY^T)(B^{-1})^T.
\]
Since the characteristic polynomials of $A$ and $(B^{-1})^T$ are relatively
prime, we have $XY^T=0$. Thus $\genfrac{[}{]}{0pt}{1} X Y \in S_{(k+l)\times n}$, which, combined with \eqref{2.26}, implies that 
$\genfrac{[}{]}{0pt}{1} X Y \in{\rm Fix}(A\oplus B, P)$.
\end{proof}

\begin{cor}\label{C2.9}
Let $f_1,\dots, f_t\in\Bbb F_2[x]\setminus\{x\}$ be irreducible such that
$\{f_1,f_1^*\},\dots$, $\{f_t,f_t^*\}$ are pairwise disjoint, where $f_i^*$ is the reciprocal polynomial of $f_i$. Let $A=A_1\oplus\cdots\oplus A_t
\in{\rm GL}(k,\Bbb F_2)$, where $A_i\in{\rm GL}(k_i,\Bbb F_2)$ whose
elementary divisors are powers of $f_i$ or $f_i^*$. Then for each $P\in\frak S_n$,
\[
|{\rm Fix}(A,P)|=\prod_{i=1}^t|{\rm Fix}(A_i,P)|.
\]
\end{cor}

\begin{lem}\label{L2.10}
Let $f\in\Bbb F_2[x]\setminus\{x\}$ be an irreducible polynomial which
is not self-reciprocal. Let $t$ be the smallest positive integer such that $f\mid x^t-1$. Let $A\in{\rm GL}(k,\Bbb F_2)$ have elementary divisors
$\underbrace{f^1,\dots,f^1}_{\mu_1},\dots,
\underbrace{f^s,\dots,f^s}_{\mu_s}$ and let $\lambda=(\lambda_1,\lambda_2,\dots)$ $\vdash n$. Then
\begin{equation}\label{2.27}
|{\rm Fix}(A,P_\lambda)|=|\{X\in M_{k\times n}: AX=XP_\lambda\}|
=2^{\deg f\sum_{j\ge 1}\lambda_{jt}\sum_{l\ge 1}
\mu_l\min\{l,2^{\nu(j)}\}}.
\end{equation}
\end{lem}

\begin{proof}
Since the characteristic polynomials of $A$ and $A^{-1}$ are relatively prime,
by the proof of Lemma~\ref{L2.8}, $AX=XP_\lambda$ implies $XX^T=0$. Hence
\[
{\rm Fix}(A,P_\lambda)=\{X\in M_{k\times n}: AX=XP_\lambda\}.
\]
To see the second equality in \eqref{2.27}, note from \eqref{2.11} that
\[
\dim \{X\in M_{k\times n}: AX=XP_\lambda\}=\sum_{i\ge 1}\bigl[
k-{\rm rank}\,(A^i-I)\bigr]\lambda_i.
\]
By \cite[Lemma 5.2 and its proof]{Hou05}, we have
\[
k-{\rm rank}\,(A^i-I)=
\begin{cases}
\displaystyle \deg f \sum_{l\ge 1}\mu_l\min\{l,2^{\nu(i)}\} &\text{if}\ t\mid i,\vspace{2pt}\cr
0 &\text{if}\ t\nmid i.\cr
\end{cases}
\]
Therefore,
\[
\dim \{X\in M_{k\times n}: AX=XP_\lambda\}=
\deg f\sum_{j\ge 1}\lambda_{jt}\sum_{l\ge 1}
\mu_l\min\{l,2^{\nu(j)}\}.
\]
\end{proof}

\section{Computation of $\Psi_{\le 3,n}$}

Recall that $\mathcal C({\rm GL}(k,\Bbb F_2))$ is a set of representatives of the conjugacy classes of ${\rm GL}(k,\Bbb F_2)$. We let 
$\mathcal C({\rm GL}(3,\Bbb F_2))=\{A_1,\dots, A_6\}$ where $A_1,\dots, A_6$ are given in 
Table~\ref{Tb1}.
\begin{table}
\caption{Information about $\mathcal C({\rm GL}(3,\Bbb F_2))$}\label{Tb1}
\vskip-5mm
\[
\begin{tabular}{l|l|c}
\hline
representative & elementary divisors& $|{\rm cent}_{{\rm GL}(3,\Bbb F_2)}(\ )|$\\ \hline
$A_1=I_3$ & $x+1,\ x+1,\ x+1$ & $2^3\cdot 3\cdot 7$ \\ \hline
$A_2=[1]\oplus\left[
\begin{matrix}
0&1\cr
1&0\cr
\end{matrix}\right]$ & $x+1,\ (x+1)^2$ & $2^3$ \\ \hline
$A_3=\left[
\begin{matrix}
0&1&0\cr
0&0&1\cr
1&1&1\cr
\end{matrix}\right]$ & $(x+1)^3$ & $2^2$ \\ \hline
$A_4=[1]\oplus\left[
\begin{matrix}
0&1\cr
1&1\cr
\end{matrix}\right]$ & $x+1,\ x^2+x+1$ & $3$ \\ \hline
$A_5=\left[
\begin{matrix}
0&1&0\cr
0&0&1\cr
1&1&0\cr
\end{matrix}\right]$ & $x^3+x+1$ & $7$ \\ \hline
$A_6=\left[
\begin{matrix}
0&1&0\cr
0&0&1\cr
1&0&1\cr
\end{matrix}\right]$ & $x^3+x^2+1$ & $7$ \\ \hline
\end{tabular}
\]
\end{table}
By \eqref{2.7}, to determine $\Psi_{\le 3,n}$, it suffices to determine $|{\rm Fix}(A_i,P_\lambda)|$ for $i=1,\dots, 6$ and for all $\lambda=(\lambda_1,\lambda_2,\dots)\vdash n$. In the following computations, we will use the notation of Theorem~\ref{T2.3} and Algorithm~\ref{A1}. 


\subsection{Computation of $|{\rm Fix}(A_1,P_\lambda)|$}\

By Lemma~\ref{L2.4} and \eqref{2.24.0},
\[
\begin{split}
|{\rm Fix}(A_1,P_\lambda)|\,&=2^{3\lambda_{0,2}}\bigl| S_{3\times \lambda_{1,2}}\bigr|\cr
&=2^{3\lambda_{0,2}}\bigl[2^{3\lambda_{1,2}-6}\bigl(36-35\delta(\lambda_{1,2})\bigr)+2^{2\lambda_{1,2}-6}\, 7\bigl(3+(-1)^{\lambda_{1,2}}\bigr)\bigr].
\end{split}
\]


\subsection{Computation of $|{\rm Fix}(A_2,P_\lambda)|$}\

We have $o(A_2)=2$. Since 
\[
A_2-I=[0]\oplus
\left[
\begin{matrix}
1&1\cr
1&1\cr
\end{matrix}\right],\qquad\text{and}\qquad 
A_2^2-I=0,
\]
we have $s_1=2$, $s_2=3$, and we can choose
\[
B_1=[1]\oplus\left[
\begin{matrix}
1\cr
1\cr
\end{matrix}\right],\quad B_2=I_3.
\]
By \eqref{2.8} and \eqref{I1.1},
\begin{equation}\label{4.2}
\alpha_1=\sum_{\nu(i)=0}\lambda_i=\lambda_{1,2},\quad
\alpha_2=\sum_{\nu(i)=1}\lambda_i=\lambda_{2,4},\quad \theta=2\alpha_1+3\alpha_2.
\end{equation}
With $Y_dY_d^T=[z_{ij}^{(d)}]$, equation \eqref{2.10} becomes
\[
\left[
\begin{matrix}
z_{11}^{(1)}&z_{12}^{(1)}&z_{12}^{(1)}\cr
z_{12}^{(1)}&z_{22}^{(1)}&z_{22}^{(1)}\cr
z_{12}^{(1)}&z_{22}^{(1)}&z_{22}^{(1)}
\end{matrix}\right]+\left[
\begin{matrix}
0& z_{12}^{(2)}+z_{13}^{(2)}& z_{12}^{(2)}+z_{13}^{(2)}\cr
 z_{12}^{(2)}+z_{13}^{(2)}&z_{22}^{(2)}+z_{33}^{(2)}&0\cr
 z_{12}^{(2)}+z_{13}^{(2)}& 0 &z_{22}^{(2)}+z_{33}^{(2)}\cr
\end{matrix}\right]=0,
\]
i.e.,
\[
\begin{cases}
z_{11}^{(1)}=z_{22}^{(1)}=0,\cr
z_{12}^{(1)}+z_{12}^{(2)}+z_{13}^{(2)}=0,\cr
z_{22}^{(2)}+z_{33}^{(2)}=0.
\end{cases}
\]
Thus, $L=4$ and
\begin{align*}
&C_1^{(1)}=\left[
\begin{matrix}
1&0\cr
0&0
\end{matrix}\right],& &C_1^{(2)}=0,\cr
&C_2^{(1)}=\left[
\begin{matrix}
0&0\cr
0&1
\end{matrix}\right],& &C_2^{(2)}=0,\cr
&C_3^{(1)}=\left[
\begin{matrix}
0&1\cr
0&0
\end{matrix}\right],& &C_3^{(2)}=\left[
\begin{matrix}
0&1&1\cr
0&0&0\cr
0&0&0
\end{matrix}\right],\cr
&C_4^{(1)}=0,& &C_4^{(2)}=\left[
\begin{matrix}
0&0&0\cr
0&1&0\cr
0&0&1
\end{matrix}\right].
\end{align*}
We have
\[
\text{type}(a_1C_1^{(1)}+\cdots+a_4C_4^{(1)})=
\begin{cases}
(2,0,0,0)&\text{if}\ (a_1,a_2,a_3)=(0,0,0),\cr
(2,0,0,1)&\text{if}\ a_3=0,\ (a_1,a_2)\ne(0,0),\cr
(2,1,0,0)&\text{if}\ a_3=1,\ (a_1,a_2)\ne(1,1),\cr
(2,1,1,0)&\text{if}\ (a_1,a_2,a_3)=(1,1,1),
\end{cases}
\]
\[
\text{type}(a_1C_1^{(2)}+\cdots+a_4C_4^{(2)})=
\begin{cases}
(3,0,0,0)&\text{if}\ (a_3,a_4)=(0,0),\cr
(3,0,0,1)&\text{if}\ (a_3,a_4)=(0,1),\cr
(3,1,0,0)&\text{if}\ a_3=1.
\end{cases}
\]
As $(a_1,\dots,a_4)$ runs over $\Bbb F_2^4$, 
\[
\begin{split}
&\text{type}(\frak C_{(a_1,\dots,a_4)})\cr
=\,&\alpha_1*\text{type}(a_1C_1^{(1)}+\cdots+a_4C_4^{(1)})
\boxplus \alpha_2*\text{type}(a_1C_1^{(2)}+\cdots+a_4C_4^{(2)})\cr
=\,&
\begin{cases}
(\theta,0,0,0)&1\ \text{time},\cr
(\theta,0,0,\delta(\alpha_2))&1\ \text{time},\cr
(\theta,0,0,\delta(\alpha_1))&3\ \text{times},\cr
(\theta,0,0,\delta(\alpha_1+\alpha_2))&3\ \text{times},\cr
(\theta,\alpha_1+\alpha_2,0,0) &6\ \text{times},\cr
(\theta,\alpha_1+\alpha_2,\alpha_1,0) &2\ \text{times}.
\end{cases}
\end{split}
\]
Therefore,
\[
\begin{split}
\frak n(A_2)
=\,&-2^\theta+2^{-3}\Bigl[2^\theta+2^{\theta-1}+(1-\delta(\alpha_2))2^{\theta-1}
+3\bigl[ 2^{\theta-1}+(1-\delta(\alpha_1))2^{\theta-1}\bigr]\cr
&+3\bigl[ 2^{\theta-1}+(1-\delta(\alpha_1+\alpha_2))2^{\theta-1}\bigr]
+6\bigl(2^{\theta-1}+2^{\theta-1-(\alpha_1+\alpha_2)}\bigr)\cr
&+2\bigl[2^{\theta-1}+(-1)^{\alpha_1}2^{\theta-1-(\alpha_1+\alpha_2)}\bigr]
\Bigr]\cr
\hfil=\,&2^{\theta-4}\bigl[ 8-3\delta(\alpha_1)-\delta(\alpha_2)-3\delta(\alpha_1+\alpha_2)\bigr] +2^{\alpha_1+2\alpha_2-3}(3+(-1)^{\alpha_1}).
\end{split}
\]
By \eqref{2.9} and \eqref{4.2},
\[
\begin{split}
&|\text{Fix}(A_2,P_\lambda)|\cr
=\,&2^{3\sum_{\nu(i)>1}\lambda_i}\, \frak n(A_2)\cr
=\,&2^{3\lambda_{0,4}}
\Bigl[ 2^{2\lambda_{1,2}+3\lambda_{2,4}-4}\bigl[ 8-3\delta(\lambda_{1,2})-\delta(\lambda_{2,4})-3\delta(\lambda_{1,2}+\lambda_{2,4})\bigr]\cr
&+ 2^{\lambda_{1,2}+2\lambda_{2,4}-3} (3+(-1)^{\lambda_{1,2}})\Bigr].
\end{split}
\]


\subsection{Computation of $|{\rm Fix}(A_3,P_\lambda)|$}\

We have $o(A_3)=4$. Since
\[
A_3-I=
\left[
\begin{matrix}
1&1&0\cr
0&1&1\cr
1&1&0\cr
\end{matrix}\right],\qquad
A_3^2-I=
\left[
\begin{matrix}
1&0&1\cr
1&0&1\cr
1&0&1\cr
\end{matrix}\right],\qquad A_3^4-I=0,
\]
we have $s_1=1, s_2=2, s_4=3$ and we can choose
\[
B_1=\left[
\begin{matrix}
1\cr 1\cr 1\cr
\end{matrix}\right],\quad
B_2=\left[
\begin{matrix}
1&0\cr 0&1\cr 1&0\cr
\end{matrix}\right],\quad B_4=I_3.
\]
By \eqref{2.8} and \eqref{I1.1}
\begin{equation}\label{4.8}
\alpha_1=\lambda_{1,2},\quad
\alpha_2=\lambda_{2,4},\quad
\alpha_4=\lambda_{4,8},\quad \theta=\alpha_1+2\alpha_2+3\alpha_4.
\end{equation}
With $Y_dY_d^T=[z_{ij}^{(d)}]$, equation \eqref{2.10} becomes
\[
\begin{split}
&\left[
\begin{matrix}
z_{11}^{(1)}&z_{11}^{(1)}&z_{11}^{(1)}\cr
z_{11}^{(1)}&z_{11}^{(1)}&z_{11}^{(1)}\cr
z_{11}^{(1)}&z_{11}^{(1)}&z_{11}^{(1)}\cr
\end{matrix}\right]+\left[
\begin{matrix}
z_{11}^{(2)}+z_{22}^{(2)}&0&z_{11}^{(2)}+z_{22}^{(2)}\cr
0&z_{11}^{(2)}+z_{22}^{(2)}&0\cr
z_{11}^{(2)}+z_{22}^{(2)}&0&z_{11}^{(2)}+z_{22}^{(2)}
\end{matrix}\right]\cr
\cr
&+\left[
\begin{matrix}
0&z_{11}^{(4)}+z_{33}^{(4)} &0\cr
z_{11}^{(4)}+z_{33}^{(4)}&0&z_{11}^{(4)}+z_{33}^{(4)}\cr
0&z_{11}^{(4)}+z_{33}^{(4)} &0
\end{matrix}\right]=0,
\end{split}
\]
i.e.,
\[
\begin{cases}
z_{11}^{(1)}+z_{11}^{(2)}+z_{22}^{(2)}=0,\cr
z_{11}^{(1)}+z_{11}^{(4)}+z_{33}^{(4)}=0.
\end{cases}
\]
Thus, $L=2$ and
\begin{align*}
&C_1^{(1)}=[1],& &C_1^{(2)}=\left[\begin{matrix} 1&0\cr 0&1\end{matrix}\right],& &C_1^{(4)}=0,\cr
&C_2^{(1)}=[1],& &C_2^{(2)}=0,& &C_2^{(4)}=\left[
\begin{matrix}
1&0&0\cr
0&0&0\cr
0&0&1
\end{matrix}\right].
\end{align*}
We have
\[
\text{type}(a_1C_1^{(1)}+a_2C_2^{(1)})=
\begin{cases}
(1,0,0,0)&\text{if}\ a_1=a_2,\cr
(1,0,0,1)&\text{if}\ a_1\ne a_2,
\end{cases}
\]
\[
\text{type}(a_1C_1^{(2)}+a_2C_2^{(2)})=
\begin{cases}
(2,0,0,0)&\text{if}\ a_1=0,\cr
(2,0,0,1)&\text{if}\ a_1=1,
\end{cases}
\]
\[
\text{type}(a_1C_1^{(4)}+a_2C_2^{(4)})=
\begin{cases}
(3,0,0,0)&\text{if}\ a_2=0,\cr
(3,0,0,1)&\text{if}\ a_2=1,
\end{cases}
\]
\[
\text{type}(\frak C_{(a_1,a_2)})=
\begin{cases}
(\theta,0,0,0)&\text{if}\ (a_1,a_2)=(0,0),\cr
(\theta,0,0,\delta(\alpha_2+\alpha_4))&\text{if}\ (a_1,a_2)=(1,1),\cr
(\theta,0,0,\delta(\alpha_1+\alpha_2))&\text{if}\ (a_1,a_2)=(1,0),\cr
(\theta,0,0,\delta(\alpha_1+\alpha_4))&\text{if}\ (a_1,a_2)=(0,1).
\end{cases}
\]
Thus,
\[
\begin{split}
\frak n(A_3)=\,&-2^\theta+2^{-1}\Bigl[2^\theta+2^{\theta-1}+\bigl(1-\delta(\alpha_2+\alpha_4)\bigr)2^{\theta-1}\cr
&+2^{\theta-1}+\bigl(1-\delta(\alpha_1+\alpha_2)\bigr)2^{\theta-1}
+2^{\theta-1}+\bigl(1-\delta(\alpha_1+\alpha_4)\bigr)2^{\theta-1}\Bigr]\cr
\hfil=\,& 2^{\theta-2} \bigl[ 4-\delta(\alpha_1+\alpha_2)-\delta(\alpha_1+\alpha_4)-\delta(\alpha_2+\alpha_4)\bigr].
\end{split}
\]
By \eqref{2.9} and \eqref{4.8},
\[
\begin{split}
&|\text{Fix}(A_3,P_\lambda)|\cr
=\,&2^{3\sum_{\nu(i)>2}\lambda_i}\,\frak n(A_3)\cr
=\,&2^{\lambda_{1,2}+2\lambda_{2,4}+3\lambda_{0,4}-2}\bigl[ 4-\delta(\lambda_{1,2}+\lambda_{2,4})-\delta(\lambda_{1,2}+\lambda_{4,8})
-\delta(\lambda_{2,4}+\lambda_{4,8})\bigr].
\end{split}
\]


\subsection{Computation of $|{\rm Fix}(A_4,P_\lambda)|$}\

By Corollary~\ref{C2.9},
\[
|{\rm Fix}(A_4,P_\lambda)|=|{\rm Fix}([1],P_\lambda)|\, 
\bigl|{\rm Fix}(\genfrac{[}{]}{0pt}{1}{0\ 1}{1\ 1}, P_\lambda) \bigr|,
\]
where $\bigl|{\rm Fix}(\genfrac{[}{]}{0pt}{1}{0\ 1}{1\ 1}, P_\lambda) \bigr|$ is given by \eqref{I5}.
By Lemma~\ref{L2.4},
\begin{equation}\label{4.10.2}
|{\rm Fix}([1],P_\lambda)|=2^{\lambda_{0,2}} \bigl|S_{1\times \lambda_{1,2}}\bigr|=
2^{\lambda_{0,1}-1}\bigl(2-\delta(\lambda_{1,2})\bigr).
\end{equation}
Therefore,
\[
|{\rm Fix}(A_4,P_\lambda)|= 2^{\lambda_{0,1}+2\lambda_{0,6}-2}\bigl(2-\delta(\lambda_{1,2})\bigr)\Bigl[2^{2\lambda_{3,6}}
+(-1)^{\lambda_{3,6}} 2^{\lambda_{3,6}}\Bigr].
\]


\subsection{Computation of $|{\rm Fix}(A_5,P_\lambda)|$ and $|{\rm Fix}(A_6,P_\lambda)|$}\

The smallest positive integer $t$ such that $x^3+x+1\mid x^t-1$ is 7. By Lemma~\ref{L2.10},
\begin{equation}\label{4.13}
|{\rm Fix}(A_5,P_\lambda)|=2^{3\lambda_{0,7}}.
\end{equation}
In the same way,
\[
|{\rm Fix}(A_6,P_\lambda)|=2^{3\lambda_{0,7}}.
\]


\section{Computation of $\Psi_{\le 4, n}$}

We let $\mathcal C\bigl({\rm GL}(4,\Bbb F_2)\bigr)=\{A_1,\dots,A_{14}\}$ where $A_1,\dots,A_{14}$ are given in Table 2.
It suffices to determine $|{\rm Fix}(A_i,P_\lambda)|$ for $i=1,\dots,14$ and for all $\lambda=(\lambda_1,\lambda_2,
\dots)\vdash n$. Again, in the following computations, we will use the notation of Theorem~\ref{T2.3} and Algorithm~\ref{A1}.

\begin{table}\label{Tb2}
\caption{Information about $\mathcal C({\rm GL}(4,\Bbb F_2))$}
\vskip-5mm
\[
\begin{tabular}{l|l|c}
\hline
representative & elementary divisors& $|{\rm cent}_{{\rm GL}(4,\Bbb F_2)}(\ )|$\\ \hline
$A_1=I_4$ & $x+1,\ x+1,\ x+1,\ x+1$ & $2^6\cdot 3^2\cdot 5\cdot 7$ \\ \hline
$A_2=I_2\oplus\left[
\begin{matrix}
0&1\cr
1&0\cr
\end{matrix}\right]$ & $x+1,\ x+1,\ (x+1)^2$ & $2^6\cdot 3$ \\ \hline
$A_3=[1]\oplus\left[
\begin{matrix}
0&1&0\cr
0&0&1\cr
1&1&1\cr
\end{matrix}\right]$ & $x+1,\ (x+1)^3$ & $2^4$ \\ \hline
$A_4=I_2\oplus\left[
\begin{matrix}
0&1\cr
1&1\cr
\end{matrix}\right]$ & $x+1,\ x+1,\ x^2+x+1$ & $2\cdot 3^2$ \\ \hline
$A_5=[1]\oplus\left[
\begin{matrix}
0&1&0\cr
0&0&1\cr
1&1&0\cr
\end{matrix}\right]$ & $x+1,\ x^3+x+1$ & $7$ \\ \hline
$A_6=[1]\oplus\left[
\begin{matrix}
0&1&0\cr
0&0&1\cr
1&0&1\cr
\end{matrix}\right]$ & $x+1,\ x^3+x^2+1$ & $7$ \\ \hline
$A_7=\left[
\begin{matrix}
0&1\cr
1&0\cr
\end{matrix}\right]\oplus\left[
\begin{matrix}
0&1\cr
1&0\cr
\end{matrix}\right]$ & $(x+1)^2,\ (x+1)^2$ & $2^5\cdot 3$ \\ \hline
$A_8=\left[
\begin{matrix}
0&1\cr
1&0\cr
\end{matrix}\right]\oplus\left[
\begin{matrix}
0&1\cr
1&1\cr
\end{matrix}\right]$ & $(x+1)^2,\ x^2+x+1$ & $2\cdot 3$ \\ \hline
$A_9=\left[
\begin{matrix}
0&1\cr
1&1\cr
\end{matrix}\right]\oplus\left[
\begin{matrix}
0&1\cr
1&1\cr
\end{matrix}\right]$ & $x^2+x+1,\ x^2+x+1$ & $2^2\cdot 3^2\cdot 5$ \\ \hline
$A_{10}=\left[
\begin{matrix}
0&1&0&0\cr
0&0&1&0\cr
0&0&0&1\cr
1&0&0&0\cr
\end{matrix}\right]$ & $(x+1)^4$ & $2^3$ \\ \hline
$A_{11}=\left[
\begin{matrix}
0&1&0&0\cr
0&0&1&0\cr
0&0&0&1\cr
1&0&1&0\cr
\end{matrix}\right]$ & $(x^2+x+1)^2$ & $2^2\cdot 3$ \\ \hline
$A_{12}=\left[
\begin{matrix}
0&1&0&0\cr
0&0&1&0\cr
0&0&0&1\cr
1&1&1&1\cr
\end{matrix}\right]$ & $x^4+x^3+x^2+x+1$ & $3\cdot 5$ \\ \hline
$A_{13}=\left[
\begin{matrix}
0&1&0&0\cr
0&0&1&0\cr
0&0&0&1\cr
1&1&0&0\cr
\end{matrix}\right]$ & $x^4+x+1$ & $3\cdot 5$ \\ \hline
$A_{14}=\left[
\begin{matrix}
0&1&0&0\cr
0&0&1&0\cr
0&0&0&1\cr
1&0&0&1\cr
\end{matrix}\right]$ & $x^4+x^3+1$ & $3\cdot 5$ \\ \hline
\end{tabular}
\]
\vskip1cm
\end{table}

\subsection{Computation of $|{\rm Fix}(A_1,P_\lambda)|$}\

By Lemma~\ref{L2.4} and \eqref{2.24.0},
\begin{equation}\label{5.1}
\begin{split}
|{\rm Fix}(A_1,P_\lambda)|=\,&2^{4\lambda_{0,2}}\bigl| S_{4\times \lambda_{1,2}} \bigr|\cr
=\,&2^{4\lambda_{0,2}}\Bigl[ 2^{4\lambda_{1,2}-10}\bigl(436-435\delta(\lambda_{1,2})\bigr)\cr
&+2^{3\lambda_{1,2}-10}\, 35\bigl(3+(-1)^{\lambda_{1,2}}\bigr)+2^{2\lambda_{1,2}-7}\, 7\bigl(5+3(-1)^{\lambda_{1,2}}\bigr)\Bigr].
\end{split}
\end{equation}


\subsection{Computation of $|{\rm Fix}(A_2,P_\lambda)|$}\

We have $o(A_2)=2$ and
\[
A_2-I=
\left[
\begin{matrix}
0&0\cr
0&0\end{matrix}\right]\oplus\left[
\begin{matrix}
1&1\cr
1&1
\end{matrix}\right],\quad 
A_2^2-I=0.
\]
Hence $s_1=3,\ s_2=4$. We can choose
\[
B_1=I_2\oplus\left[
\begin{matrix}
1\cr
1\cr
\end{matrix}\right],\quad B_2=I_4.
\]
By \eqref{2.8} and \eqref{I1.1},
\begin{equation}\label{5.2}
\alpha_1=\lambda_{1,2},\quad
\alpha_2=\lambda_{2,4},\quad \theta=3\alpha_1+4\alpha_2.
\end{equation}
With $Y_dY_d^T=[z_{ij}^{(d)}]$, equation \eqref{2.10} becomes
\[
\left[
\begin{matrix}
z_{11}^{(1)}&z_{12}^{(1)}&z_{13}^{(1)}&z_{13}^{(1)}\cr
z_{12}^{(1)}&z_{22}^{(1)}&z_{23}^{(1)}&z_{23}^{(1)}\cr
z_{13}^{(1)}&z_{23}^{(1)}&z_{33}^{(1)}&z_{33}^{(1)}\cr
z_{13}^{(1)}&z_{23}^{(1)}&z_{33}^{(1)}&z_{33}^{(1)}\cr
\end{matrix}\right]+ \left[
\begin{matrix}
0&0& z_{13}^{(2)}+z_{14}^{(2)}& z_{13}^{(2)}+z_{14}^{(2)}\cr
0&0& z_{23}^{(2)}+z_{24}^{(2)}& z_{23}^{(2)}+z_{24}^{(2)}\cr
z_{13}^{(2)}+z_{14}^{(2)}& z_{23}^{(2)}+z_{24}^{(2)}&z_{33}^{(2)}+z_{44}^{(2)} &0\cr
z_{13}^{(2)}+z_{14}^{(2)}& z_{23}^{(2)}+z_{24}^{(2)}&0&z_{33}^{(2)}+z_{44}^{(2)}\cr
\end{matrix}\right]=0,
\]
i.e., 
\[
\begin{cases}
z_{11}^{(1)}=z_{22}^{(1)}=z_{33}^{(1)}=z_{12}^{(1)}=0,\cr
z_{13}^{(1)}+z_{13}^{(2)}+z_{14}^{(2)}=0,\cr
z_{23}^{(1)}+z_{23}^{(2)}+z_{24}^{(2)}=0,\cr
z_{33}^{(2)}+z_{44}^{(2)}=0.
\end{cases}
\]
Thus, $L=7$ and
\begin{align*}
&C_1^{(1)}=\left[
\begin{matrix}
1&0&0\cr
0&0&0\cr
0&0&0\end{matrix}\right],& &C_1^{(2)}=0,\cr
&C_2^{(1)}=\left[
\begin{matrix}
0&0&0\cr
0&1&0\cr
0&0&0\end{matrix}\right],& &C_2^{(2)}=0,\cr
&C_3^{(1)}=\left[
\begin{matrix}
0&0&0\cr
0&0&0\cr
0&0&1\end{matrix}\right],& &C_3^{(2)}=0,\cr
&C_4^{(1)}=\left[
\begin{matrix}
0&1&0\cr
0&0&0\cr
0&0&0\end{matrix}\right],& &C_4^{(2)}=0,\cr
&C_5^{(1)}=\left[
\begin{matrix}
0&0&1\cr
0&0&0\cr
0&0&0\end{matrix}\right],& &C_5^{(2)}=\left[
\begin{matrix}
0&0&1&1\cr
0&0&0&0\cr
0&0&0&0\cr
0&0&0&0
\end{matrix}\right],\cr
&C_6^{(1)}=\left[
\begin{matrix}
0&0&0\cr
0&0&1\cr
0&0&0\end{matrix}\right],& &C_6^{(2)}=\left[
\begin{matrix}
0&0&0&0\cr
0&0&1&1\cr
0&0&0&0\cr
0&0&0&0
\end{matrix}\right],\cr
&C_7^{(1)}=0,& &C_7^{(2)}=\left[
\begin{matrix}
0&0&0&0\cr
0&0&0&0\cr
0&0&1&0\cr
0&0&0&1
\end{matrix}\right].
\end{align*}
Fix $(a_4,a_5,a_6,a_7)\in\Bbb F_2^4$ and let $(a_1,a_2,a_3)$ run over $\Bbb F_2^3$.
If $(a_4,a_5,a_6)=(0,0,0)$, 
\[
\text{type}(a_1C_1^{(1)}+\cdots+a_7C_7^{(1)})=
\begin{cases}
(3,0,0,0)& 1\ \text{time},\cr
(3,0,0,1)& 7\ \text{times}.
\end{cases}
\]
If $(a_4,a_5,a_6)\ne(0,0,0)$, 
\[
\text{type}(a_1C_1^{(1)}+\cdots+a_7C_7^{(1)})=
\begin{cases}
(3,1,0,0)& 3\ \text{times},\cr
(3,1,1,0)& 1\ \text{time},\cr
(3,\square, 0, 1)& 4\ \text{times}.
\end{cases}
\]
In the above, the symbol $\square$ represents a component of a type which is not needed in the subsequent computations. We also have
\[
\text{type}(a_1C_1^{(2)}+\cdots+a_7C_7^{(2)})=
\begin{cases}
(4,0,0,0)&\text{if}\ (a_5,a_6,a_7)=(0,0,0),\cr
(4,0,0,1)&\text{if}\ (a_5,a_6,a_7)=(0,0,1),\cr 
(4,1,0,0)&\text{if}\ (a_5,a_6)\ne (0,0).
\end{cases}
\]
Therefore, as $(a_1,\dots,a_7)$ runs over $\Bbb F_2^7$,
\[
\begin{split}
&\text{type}(\frak C_{(a_1,\dots,a_7)})\cr
=\,&\alpha_1*\text{type}(a_1C_1^{(1)}+\cdots+a_7C_7^{(1)})\boxplus
\alpha_2*\text{type}(a_1C_1^{(2)}+\cdots+a_7C_7^{(2)})\cr
=\,&
\begin{cases}
(\theta,0,0,0)& 1\ \text{time},\cr
(\theta,0,0,\delta(\alpha_2))& 1\ \text{time},\cr
(\theta,\alpha_1\square,0,\delta(\alpha_1))& 11\ \text{times},\cr
(\theta,\alpha_1\square,0,\delta(\alpha_1+\alpha_2))& 11\ \text{times},\cr
(\theta,\alpha_1,0,0)& 3\ \text{times},\cr
(\theta,\alpha_1,0,\delta(\alpha_2))& 3\ \text{times},\cr
(\theta,\alpha_1+\alpha_2,0,0)& 2^2\cdot3^2\ \text{times},\cr
(\theta,\alpha_1,\alpha_1,0)& 1\ \text{time},\cr
(\theta,\alpha_1,\alpha_1,\delta(\alpha_2))& 1\ \text{time},\cr
(\theta,\alpha_1+\alpha_2,\alpha_1,0)& 2^2\cdot 3\ \text{times},\cr
(\theta,\alpha_1\square+\alpha_2,0,\delta(\alpha_1))& 2^4\cdot 3\ \text{times}.
\end{cases}
\end{split}
\]
By \eqref{n(A)},
\[
\begin{split}
\frak n(A_2)=&\;2^{3\alpha_1+4\alpha_2-7}\bigl[ 24-\delta(\alpha_2)-11\delta(\alpha_1)-11\delta(\alpha_1+\alpha_2)\bigr]\cr
&+2^{2\alpha_1+4\alpha_2-7}\bigl(3+(-1)^{\alpha_1}\bigr)\bigl(2-\delta(\alpha_2)\bigr)\cr
&+2^{2\alpha_1+3\alpha_2-5}\bigl[7+(-1)^{\alpha_1}-4\delta(\alpha_1)\bigr].
\end{split}
\]
By \eqref{2.9} and \eqref{5.2},
\[
\begin{split}
&|\text{Fix}(A_2,P_\lambda)|\cr
=\,&2^{4 \sum_{\nu(i)>1}\lambda_i}\,\frak n(A_2)\cr
=\,&2^{4\lambda_{0,4}}\Bigl[2^{3\lambda_{1,2}+4\lambda_{2,4}-7}\bigl[ 24-\delta(\lambda_{2,4})-11\delta(\lambda_{1,2})-11\delta(\lambda_{1,2}
+\lambda_{2,4})\bigr]\cr
&\!+2^{2\lambda_{1,2}+4\lambda_{2,4}-7}\bigl(3+(-1)^{\lambda_{1,2}}\bigr)\bigl(2-\delta(\lambda_{2,4})\bigr)\cr
&+2^{2\lambda_{1,2}+3\lambda_{2,4}-5}\, 3 \bigl[ 7+(-1)^{\lambda_{1,2}}-4\delta(\lambda_{1,2})\bigr]\Bigr].
\end{split}
\]


\subsection{Computation of $|{\rm Fix}(A_3,P_\lambda)|$}\

The data for the matrix $A_3$ are as follows:
$o(A_3)=4$, $s_1=2$, $s_2=3$, $s_4=4$,
\[
B_1=[1]\oplus\left[
\begin{matrix}
1\cr
1\cr
1\cr
\end{matrix}\right],\quad 
B_2=[1]\oplus\left[
\begin{matrix}
1&0\cr
0&1\cr
1&0\cr
\end{matrix}\right],\quad B_4=I_4,
\]
\begin{equation}\label{5.2.1}
\alpha_1=\lambda_{1,2},\quad
\alpha_2=\lambda_{2,4},\quad
\alpha_4=\lambda_{4,8},\quad \theta=2\alpha_1+3\alpha_2+4\alpha_4.
\end{equation}
With $Y_dY_d^T=[z_{ij}^{(d)}]$, equation~\eqref{2.10} becomes
\[
\begin{split}
&\left[
\begin{matrix}
z_{11}^{(1)}&z_{12}^{(1)}&z_{12}^{(1)}&z_{12}^{(1)}\cr
z_{12}^{(1)}&z_{22}^{(1)}&z_{22}^{(1)}&z_{22}^{(1)}\cr
z_{12}^{(1)}&z_{22}^{(1)}&z_{22}^{(1)}&z_{22}^{(1)}\cr
z_{12}^{(1)}&z_{22}^{(1)}&z_{22}^{(1)}&z_{22}^{(1)}\cr
\end{matrix}\right]+ \left[
\begin{matrix}
0&z_{12}^{(2)}+z_{13}^{(2)} &z_{12}^{(2)}+z_{13}^{(2)}&z_{12}^{(2)}+z_{13}^{(2)}\cr
z_{12}^{(2)}+z_{13}^{(2)}&z_{22}^{(2)}+z_{33}^{(2)}&0&z_{22}^{(2)}+z_{33}^{(2)}\cr
z_{12}^{(2)}+z_{13}^{(2)}&0&z_{22}^{(2)}+z_{33}^{(2)}&0\cr
z_{12}^{(2)}+z_{13}^{(2)}&z_{22}^{(2)}+z_{33}^{(2)}&0&z_{22}^{(2)}+z_{33}^{(2)}\cr
\end{matrix}\right]\cr
\cr
&+ \left[
\begin{matrix}
0&0&0&0\cr
0&0&z_{22}^{(4)}+z_{44}^{(4)}&0\cr
0&z_{22}^{(4)}+z_{44}^{(4)}&0&z_{22}^{(4)}+z_{44}^{(4)}\cr
0&0&z_{22}^{(4)}+z_{44}^{(4)}&0\cr
\end{matrix}\right]=0,
\end{split}
\]
i.e.,
\[
\begin{cases}
z_{11}^{(1)}=0,\cr
z_{12}^{(1)}+z_{12}^{(2)}+z_{13}^{(2)}=0,\cr
z_{22}^{(1)}+z_{22}^{(2)}+z_{33}^{(2)}=0,\cr
z_{22}^{(1)}+z_{22}^{(4)}+z_{44}^{(4)}=0.
\end{cases}
\]
Thus, $L=4$ and
\begin{align*}
&C_1^{(1)}=\left[\begin{matrix} 1&0\cr 0&0\end{matrix}\right],&
&C_1^{(2)}=0,&
&C_1^{(4)}=0,\\
&C_2^{(1)}=\left[\begin{matrix} 0&1\cr 0&0\end{matrix}\right],&
&C_2^{(2)}=\left[
\begin{matrix}
0&1&1\cr
0&0&0\cr
0&0&0\end{matrix}\right],&
&C_2^{(4)}=0,\\
&C_3^{(1)}=\left[\begin{matrix} 0&0\cr 0&1\end{matrix}\right],&
&C_3^{(2)}=\left[
\begin{matrix}
0&0&0\cr
0&1&0\cr
0&0&1\end{matrix}\right],&
&C_3^{(4)}=0,\\
&C_4^{(1)}=\left[\begin{matrix} 0&0\cr 0&1\end{matrix}\right],&
&C_4^{(2)}=0,&
&C_3^{(4)}=\left[
\begin{matrix}
0&0&0&0\cr
0&1&0&0\cr
0&0&0&0\cr
0&0&0&1
\end{matrix}\right].
\end{align*}
We have
\[
\begin{split}
&\text{type}(a_1C_1^{(1)}+\cdots+a_4C_4^{(1)})\cr
=\,&
\begin{cases}
(2,0,0,0)&\text{if}\ a_2=0,\ (a_1,a_3,a_4)=(0,0,0), (0,1,1),\cr
(2,0,0,1)&\text{if}\ a_2=0,\ (a_1,a_3,a_4)\ne(0,0,0), (0,1,1),\cr
(2,1,1,0)&\text{if}\ a_2=1,\ (a_1,a_3,a_4)=(1,0,1), (1,1,0),\cr
(2,1,0,0)&\text{if}\ a_2=0,\ (a_1,a_3,a_4)\ne(1,0,1), (1,1,0),
\end{cases}
\end{split}
\]
\[
\text{type}(a_1C_1^{(2)}+\cdots+a_4C_4^{(2)})=
\begin{cases}
(3,0,0,0)&\text{if}\ (a_2,a_3)=(0,0),\cr
(3,0,0,1)&\text{if}\ (a_2,a_3)=(0,1),\cr
(3,1,0,0)&\text{if}\ a_2=1,
\end{cases}
\]
\[
\text{type}(a_1C_1^{(4)}+\cdots+a_4C_4^{(4)})=
\begin{cases}
(4,0,0,0)&\text{if}\ a_4=0,\cr
(4,0,0,1)&\text{if}\ a_4=1.
\end{cases}
\]
As $(a_1,\dots,a_4)$ runs over $\Bbb F_2^4$, 
\[
\text{type}(\frak C_{(a_1,\dots,a_4)})=
\begin{cases}
(\theta, 0,0,0)& 1\ \text{time},\cr
(\theta, 0,0,\delta(\alpha_2+\alpha_4))& 1\ \text{time},\cr
(\theta, 0,0,\delta(\alpha_1))& 1\ \text{time},\cr
(\theta, 0,0,\delta(\alpha_1+\alpha_4))& 2\ \text{times},\cr
(\theta, 0,0,\delta(\alpha_1+\alpha_2))& 2\ \text{times},\cr
(\theta, 0,0,\delta(\alpha_1+\alpha_2+\alpha_4))& 1\ \text{time},\cr
(\theta, \alpha_1+\alpha_2,\alpha_1,0)& 1\ \text{time},\cr
(\theta, \alpha_1+\alpha_2,\alpha_1,\delta(\alpha_4))& 1\ \text{time},\cr
(\theta, \alpha_1+\alpha_2,0,0)& 3\ \text{times},\cr
(\theta, \alpha_1+\alpha_2,0,\delta(\alpha_4))& 3\ \text{times}.
\end{cases}
\]
By \eqref{n(A)},
\[
\begin{split}
\frak n(A_3)
=\,& 2^{2\alpha_1+3\alpha_2+4\alpha_4-4}
\bigl[8-\delta(\alpha_2+\alpha_4)-\delta(\alpha_1)-2\delta(\alpha_1+\alpha_4)
-2\delta(\alpha_1+\alpha_2)\cr
&\! -\delta(\alpha_1+\alpha_2+\alpha_4)\bigr]
+2^{\alpha_1+2\alpha_2+4\alpha_4-4}\bigl(3+(-1)^{\alpha_1}\bigr)\bigl(2-\delta(\alpha_4)\bigr).
\end{split}
\]
By \eqref{2.9} and \eqref{5.2.1},
\[
\begin{split}
&|\text{Fix}(A_3,P_\lambda)|\cr
=\,& 2^{4\sum_{\nu(i)>2}\lambda_i}\, \frak n(A_3)\cr
\hfil =\,& 2^{4\lambda_{0,8}}\Bigl[ 2^{2\lambda_{1,2}+3\lambda_{2,4}+4\lambda_{4,8}-4}
\bigl[8-\delta(\lambda_{2,4}+\lambda_{4,8})-\delta(\lambda_{1,2})\cr 
&\! -2\delta(\lambda_{1,2}+\lambda_{4,8})
-2\delta(\lambda_{1,2}+\lambda_{2,4})-\delta(\lambda_{1,2}+\lambda_{2,4}+\lambda_{4,8})\bigr]\cr
&\! +2^{\lambda_{1,2}+2\lambda_{2,4}+4\lambda_{4,8}-4}\bigl(3+(-1)^{\lambda_{1,2}}\bigr)\bigl(2-\delta(\lambda_{4,8})\bigr)\Bigr].
\end{split}
\]


\subsection{Computation of $|{\rm Fix}(A_4,P_\lambda)|$}\

By Corollary~\ref{C2.9}, \eqref{I5}, \eqref{2.14} and \eqref{2.24.0}, we have
\[
\begin{split}
|{\rm Fix}(A_4,P_\lambda)|=\,&|{\rm Fix}(I_2,P_\lambda)|\,|{\rm Fix}(\genfrac{[}{]}{0pt}{1}{0\ 1}{1\ 1},P_\lambda)|\cr
=\,&2^{2\lambda_{0,2}+2\lambda_{0,6}-1}\bigl( 2^{2\lambda_{3,6}}+(-1)^{\lambda_{3,6}}
2^{\lambda_{3,6}}\bigr)\cr
&\cdot\bigl[2^{2\lambda_{1,2}-3}\bigl(4-3\delta(\lambda_{1,2})\bigr)+2^{\lambda_{1,2}-3}\bigl(3+(-1)^{\lambda_{1,2}}\bigr)\bigr].
\end{split}
\]


\subsection{Computation of $|{\rm Fix}(A_5,P_\lambda)|$ and $|{\rm Fix}(A_6,P_\lambda)|$}\

Write $A_5=[1]\oplus C$ where
\[
C=\left[
\begin{matrix}
0&1&0\cr
0&0&1\cr
1&1&0\cr
\end{matrix}\right].
\]
By Corollary~\ref{C2.9},
\[
|{\rm Fix}(A_5,P_\lambda)|=|{\rm Fix}([1],P_\lambda)|\, |{\rm Fix}(C,P_\lambda)|,
\]
where $|{\rm Fix}([1],P_\lambda)|$ is given by \eqref{4.10.2} 
and $|{\rm Fix}(C,P_\lambda)|$ is given by \eqref{4.13}. Hence we have 
\[
|{\rm Fix}(A_5,P_\lambda)|=2^{\lambda_{0,1}+3\lambda_{0,7}-1}\bigl(2-\delta(\lambda_{1,2})\bigr).
\]
In the same way,
\[
|{\rm Fix}(A_6,P_\lambda)|=2^{\lambda_{0,1}+3\lambda_{0,7}-1}\bigl(2-\delta(\lambda_{1,2})\bigr).
\]


\subsection{Computation of $|{\rm Fix}(A_7,P_\lambda)|$}\

The data for $A_7$ are as follows: $o(A_7)=2$, $s_1=2$, $s_2=4$,
\[
B_1=\left[
\begin{matrix} 1\cr 1\end{matrix}\right]\oplus\left[\begin{matrix} 1\cr 1\end{matrix}\right],
\quad B_2=I_4,
\]
\begin{equation}\label{5.2.2}
\alpha_1=\lambda_{1,2},\quad
\alpha_2=\lambda_{2,4},\quad \theta=2\alpha_1+4\alpha_2.
\end{equation}
With $Y_dY_d^T=[z_{ij}^{(d)}]$, equation \eqref{2.10} becomes
\[
\left[
\begin{matrix}
z_{11}^{(1)}&z_{11}^{(1)}&z_{12}^{(1)}&z_{12}^{(1)}\cr
z_{11}^{(1)}&z_{11}^{(1)}&z_{12}^{(1)}&z_{12}^{(1)}\cr
z_{12}^{(1)}&z_{12}^{(1)}&z_{22}^{(1)}&z_{22}^{(1)}\cr
z_{12}^{(1)}&z_{12}^{(1)}&z_{22}^{(1)}&z_{22}^{(1)}\cr
\end{matrix}\right]+\left[
\begin{matrix}
z_{11}^{(2)}+z_{22}^{(2)}&0&z_{13}^{(2)}+z_{24}^{(2)}&z_{14}^{(2)}+z_{23}^{(2)}\cr
0&z_{11}^{(2)}+z_{22}^{(2)}&z_{14}^{(2)}+z_{23}^{(2)}&z_{13}^{(2)}+z_{24}^{(2)}\cr
z_{13}^{(2)}+z_{24}^{(2)}&z_{14}^{(2)}+z_{23}^{(2)}&z_{33}^{(2)}+z_{44}^{(2)}&0\cr
z_{14}^{(2)}+z_{23}^{(2)}&z_{13}^{(2)}+z_{24}^{(2)}&0&z_{33}^{(2)}+z_{44}^{(2)}\cr
\end{matrix}\right]=0,
\]
i.e.,
\[
\begin{cases}
z_{11}^{(1)}=z_{22}^{(1)}=0,\cr
z_{11}^{(2)}+z_{22}^{(2)}=0,\cr
z_{33}^{(2)}+z_{44}^{(2)}=0,\cr
z_{12}^{(1)}+z_{14}^{(2)}+z_{23}^{(2)}=0,\cr
z_{12}^{(1)}+z_{13}^{(2)}+z_{24}^{(2)}=0.
\end{cases}
\]
Thus, $L=6$ and
\begin{align*}
&C_1^{(1)}=\left[\begin{matrix} 1&0\cr 0&0\end{matrix}\right],& &C_1^{(2)}=0,\\
&C_2^{(1)}=\left[\begin{matrix} 0&0\cr 0&1\end{matrix}\right],& &C_2^{(2)}=0,\\
&C_3^{(1)}=0,& &C_3^{(2)}=\left[
\begin{matrix}
1&0&0&0\cr
0&1&0&0\cr
0&0&0&0\cr
0&0&0&0
\end{matrix}\right],\\ 
&C_4^{(1)}=0,& &C_4^{(2)}=\left[
\begin{matrix}
0&0&0&0\cr
0&0&0&0\cr
0&0&1&0\cr
0&0&0&1
\end{matrix}\right],\\ 
&C_5^{(1)}=\left[\begin{matrix} 0&1\cr 0&0\end{matrix}\right],& &C_5^{(2)}=\left[
\begin{matrix}
0&0&0&1\cr
0&0&1&0\cr
0&0&0&0\cr
0&0&0&0
\end{matrix}\right],\\
&C_6^{(1)}=\left[\begin{matrix} 0&1\cr 0&0\end{matrix}\right],& &C_6^{(2)}=\left[
\begin{matrix}
0&0&1&0\cr
0&0&0&1\cr
0&0&0&0\cr
0&0&0&0
\end{matrix}\right].
\end{align*}
We have
\[
\text{type}(a_1C_1^{(1)}+\cdots+a_6C_6^{(1)})=
\begin{cases}
(2,0,0,0)&\text{if}\ a_5=a_6,\ (a_1,a_2)=(0,0),\cr
(2,0,0,1)&\text{if}\ a_5=a_6,\ (a_1,a_2)\ne (0,0),\cr
(2,1,0,0)&\text{if}\ a_5\ne a_6,\ (a_1,a_2)\ne (1,1),\cr
(2,1,1,0)&\text{if}\ a_5\ne a_6,\ (a_1,a_2)=(1,1),
\end{cases}
\]
\[
\text{type}(a_1C_1^{(2)}+\cdots+a_6C_6^{(2)})=
\begin{cases}
(4,0,0,0)&\text{if}\ (a_5,a_6)=(0,0),\ (a_3,a_4)=(0,0),\cr
(4,0,0,1)&\text{if}\ (a_5,a_6)=(0,0),\ (a_3,a_4)\ne(0,0),\cr
(4,2,0,0)&\text{if}\ a_5\ne a_6,\cr
(4,1,0,0)&\text{if}\ (a_5,a_6)=(1,1),\ (a_3,a_4)\ne(1,1),\cr
(4,1,1,0)&\text{if}\ (a_5,a_6)=(1,1),\ (a_3,a_4)=(1,1).
\end{cases}
\]
As $(a_1,\dots,a_6)$ runs over $\Bbb F_2^6$,
\[
\text{type}(\frak C_{(a_1,\dots,a_6)})=
\begin{cases}
(\theta,0,0,0)& 1\ \text{time},\cr
(\theta,0,0,\delta(\alpha_2))& 3\ \text{times},\cr
(\theta,\alpha_2,0,0)& 3\ \text{times},\cr
(\theta,\alpha_2,\alpha_2,0)& 1\ \text{time},\cr
(\theta,0,0,\delta(\alpha_1))& 3\ \text{times},\cr
(\theta,0,0,\delta(\alpha_1+\alpha_2))& 3^2\ \text{times},\cr
(\theta,\alpha_2,0,\delta(\alpha_1))& 3^2\ \text{times},\cr
(\theta,\alpha_2,\alpha_2,\delta(\alpha_1))& 3\ \text{times},\cr
(\theta,\alpha_1+2\alpha_2,0,0)& 2^3\cdot3\ \text{times},\cr
(\theta,\alpha_1+2\alpha_2,\alpha_1,0)& 2^3\ \text{times}.
\end{cases}
\]
By \eqref{n(A)},
\[
\begin{split}
\frak n(A_7)=\,&2^{2\alpha_1+4\alpha_2-6}\bigl[16-3\delta(\alpha_2)-3\delta(\alpha_1)
-9 \delta(\alpha_1+\alpha_2)\bigl]\cr
&\! +2^{2\alpha_1+3\alpha_2-6}\bigl(3+(-1)^{\alpha_2}\bigr)\bigl(4-3\delta(\alpha_1)\bigr)
+2^{\alpha_1+2\alpha_2-3}\bigl(3+(-1)^{\alpha_1}\bigr).
\end{split}
\]
By \eqref{2.9} and \eqref{5.2.2},
\[
\begin{split}
&|\text{Fix}(A_7,P_\lambda)|\cr
=\,&2^{4\sum_{\nu(i)>1}\lambda_i}\, \frak n(A_7)\cr
=\,&2^{4\lambda_{0,4}}\Bigl[2^{2\lambda_{1,2}+4\lambda_{2,4}-6}\bigl[16-3\delta(\lambda_{2,4})-3\delta(\lambda_{1,2})-9\delta(\lambda_{1,2}+\lambda_{2,4})\bigl]\cr
&\! +2^{2\lambda_{1,2}+3\lambda_{2,4}-6}\bigl(3+(-1)^{\lambda_{2,4}}\bigr)\bigl(4-3\delta(\lambda_{1,2})\bigr)
+2^{\lambda_{1,2}+2\lambda_{2,4}-3}\bigl(3+(-1)^{\lambda_{1,2}}\bigr)\Bigr].
\end{split}
\]


\subsection{Computation of $|{\rm Fix}(A_8,P_\lambda)|$}\
 
The result follows immediately from Corollary~\ref{C2.9} and equations \eqref{I4} and \eqref{I5}. 
We have
\[
\begin{split}
|{\rm Fix}(A_8,P_\lambda)|
=\,&|{\rm Fix}(\genfrac{[}{]}{0pt}{1}{0\ 1}{1\ 0},P_\lambda)|\, 
|{\rm Fix}(\genfrac{[}{]}{0pt}{1}{0\ 1}{1\ 1},P_\lambda)|\cr
\hfil =\,& 2^{\lambda_{1,2}+2\lambda_{0,2}+2\lambda_{0,6}-3}\cr
&\cdot\bigl[ 4-\delta(\lambda_{1,2})-\delta(\lambda_{2,4})-\delta(\lambda_{1,2}+\lambda_{2,4})\bigr]\bigl[
2^{2\lambda_{3,6}}+(-1)^{\lambda_{3,6}}2^{\lambda_{3,6}}\bigr].
\end{split}
\]


\subsection{Computation of $|{\rm Fix}(A_9,P_\lambda)|$}\

We have $o(A_9)=3$, $s_1=0$, $s_3=4$, $B_1=\emptyset$, $B_3=I_4$ and
\begin{equation}\label{5.2.3}
\alpha_3=\lambda_{3,6}.
\end{equation}
With $Y_dY_d^T=[z_{ij}^{(d)}]$, equation \eqref{2.10} becomes
\[
\left[
\begin{matrix}
0&z_{11}^{(3)}+z_{12}^{(3)}+z_{22}^{(3)}& z_{14}^{(3)}+z_{23}^{(3)}& z_{13}^{(3)}+z_{14}^{(3)}+z_{24}^{(3)}\cr
z_{11}^{(3)}+z_{12}^{(3)}+z_{22}^{(3)}& 0 & z_{13}^{(3)}+z_{23}^{(3)}+z_{24}^{(3)} & z_{14}^{(3)}+z_{23}^{(3)}\cr
z_{14}^{(3)}+z_{23}^{(3)}& z_{13}^{(3)}+z_{23}^{(3)}+z_{24}^{(3)} & 0 & z_{33}^{(3)}+z_{34}^{(3)}+z_{44}^{(3)}\cr
z_{13}^{(3)}+z_{14}^{(3)}+z_{24}^{(3)}& z_{14}^{(3)}+z_{23}^{(3)}& z_{33}^{(3)}+z_{34}^{(3)}+z_{44}^{(3)}& 0
\end{matrix}\right]=0,
\]
i.e.,
\[
\begin{cases}
z_{11}^{(3)}+z_{12}^{(3)}+z_{22}^{(3)}=0,\cr
z_{33}^{(3)}+z_{34}^{(3)}+z_{44}^{(3)}=0,\cr
z_{14}^{(3)}+z_{23}^{(3)}=0,\cr
z_{13}^{(3)}+z_{14}^{(3)}+z_{24}^{(3)}=0.
\end{cases}
\]
Thus, $L=4$ and
\[
\begin{split}
C_1^{(3)}=\left[
\begin{matrix}
1&1&0&0\cr
0&1&0&0\cr
0&0&0&0\cr
0&0&0&0
\end{matrix}\right],\quad C_2^{(3)}=\left[
\begin{matrix}
0&0&0&0\cr
0&0&0&0\cr
0&0&1&1\cr
0&0&0&1
\end{matrix}\right],\cr
C_3^{(3)}=\left[
\begin{matrix}
0&0&0&1\cr
0&0&1&0\cr
0&0&0&0\cr
0&0&0&0
\end{matrix}\right],\quad C_4^{(3)}=\left[
\begin{matrix}
0&0&1&1\cr
0&0&0&1\cr
0&0&0&0\cr
0&0&0&0
\end{matrix}\right].
\end{split}
\]
As $(a_1,\dots,a_4)$ runs over $\Bbb F_2^4$,
\[
\text{type}(a_1C_1^{(3)}+\cdots+a_4 C_4^{(3)})=
\begin{cases}
(4,0,0,0)& 1\ \text{time},\cr
(4,1,1,0)& 5\ \text{times},\cr
(4,2,0,0)& 10\ \text{times},
\end{cases}
\]
\[
\text{type}(\frak C_{(a_1,\dots,a_4)})=
\begin{cases}
(4\alpha_3,0,0,0)& 1\ \text{time},\cr
(4\alpha_3,\alpha_3,\alpha_3,0)& 5\ \text{times},\cr
(4\alpha_3,2\alpha_3,0,0)& 10\ \text{times}.
\end{cases}
\]
By \eqref{n(A)},
\[
\frak n(A_9)=2^{4\alpha_3-4}+2^{3\alpha_3-4}\,5(-1)^{\alpha_3}+2^{2\alpha_3-3}\,5.
\]
By \eqref{2.9} and \eqref{5.2.3},
\[
\begin{split}
|\text{Fix}(A_9,P_\lambda)|\,&=2^{4\sum_{i\equiv 0\,(6)}\lambda_i}\,\frak n(A_9)\cr
&=2^{4\lambda_{0,6}-3}\bigl[2^{4\lambda_{3,6}-1}+2^{3\lambda_{3,6}-1}\,5(-1)^{\lambda_{3,6}}+2^{2\lambda_{3,6}}\,5\bigr].
\end{split}
\]


\subsection{Computation of $|{\rm Fix}(A_{10},P_\lambda)|$}\

We have $o(A_{10})=4$,  $s_1=1$, $s_2=2$, $s_4=4$,
\[
B_1=\left[\begin{matrix} 1\cr 1\cr 1\cr 1\end{matrix}\right],\quad
B_2=\left[\begin{matrix} 1&0\cr 0&1\cr 1&0\cr 0&1\end{matrix}\right],\quad
B_4=I_4,
\]
\begin{equation}\label{5.2.4}
\alpha_1=\lambda_{1,2},\quad \alpha_2=\lambda_{2,4},\quad 
\alpha_4=\lambda_{4,8},\quad \theta=\alpha_1+2\alpha_2+4\alpha_4.
\end{equation}
With $Y_dT_d^T=[z_{ij}^{(d)}]$, equation \eqref{2.10} becomes
\[
\begin{split}
&\left[
\begin{matrix}
z_{11}^{(1)}&z_{11}^{(1)}&z_{11}^{(1)}&z_{11}^{(1)}\cr
z_{11}^{(1)}&z_{11}^{(1)}&z_{11}^{(1)}&z_{11}^{(1)}\cr
z_{11}^{(1)}&z_{11}^{(1)}&z_{11}^{(1)}&z_{11}^{(1)}\cr
z_{11}^{(1)}&z_{11}^{(1)}&z_{11}^{(1)}&z_{11}^{(1)}
\end{matrix}\right]+\left[
\begin{matrix}
z_{11}^{(2)}+z_{22}^{(2)} & 0 & z_{11}^{(2)}+z_{22}^{(2)} & 0\cr
0 & z_{11}^{(2)}+z_{22}^{(2)} & 0 & z_{11}^{(2)}+z_{22}^{(2)}\cr
z_{11}^{(2)}+z_{22}^{(2)} & 0 & z_{11}^{(2)}+z_{22}^{(2)} & 0\cr
0 & z_{11}^{(2)}+z_{22}^{(2)} & 0 & z_{11}^{(2)}+z_{22}^{(2)}
\end{matrix}\right]\cr
\cr
+&\left[
\begin{matrix}
z_{11}^{(4)}+z_{22}^{(4)}+z_{33}^{(4)}+z_{44}^{(4)} &
z_{12}^{(4)}+z_{23}^{(4)}+z_{34}^{(4)}+z_{14}^{(4)}\cr
z_{12}^{(4)}+z_{23}^{(4)}+z_{34}^{(4)}+z_{14}^{(4)} &
z_{11}^{(4)}+z_{22}^{(4)}+z_{33}^{(4)}+z_{44}^{(4)}\cr 
0 &
z_{12}^{(4)}+z_{23}^{(4)}+z_{34}^{(4)}+z_{14}^{(4)}\cr
z_{12}^{(4)}+z_{23}^{(4)}+z_{34}^{(4)}+z_{14}^{(4)} &
0
\end{matrix}\right.\cr
\cr
&\kern2mm\left.
\begin{matrix}
0 &
z_{12}^{(4)}+z_{23}^{(4)}+z_{34}^{(4)}+z_{14}^{(4)}\cr
z_{12}^{(4)}+z_{23}^{(4)}+z_{34}^{(4)}+z_{14}^{(4)} &
0\cr
z_{11}^{(4)}+z_{22}^{(4)}+z_{33}^{(4)}+z_{44}^{(4)} &
z_{12}^{(4)}+z_{23}^{(4)}+z_{34}^{(4)}+z_{14}^{(4)}\cr
z_{12}^{(4)}+z_{23}^{(4)}+z_{34}^{(4)}+z_{14}^{(4)} &
z_{11}^{(4)}+z_{22}^{(4)}+z_{33}^{(4)}+z_{44}^{(4)}
\end{matrix}\right]=0,
\end{split}
\]
i.e.,
\[
\begin{cases}
z_{11}^{(1)}+z_{11}^{(2)}+z_{22}^{(2)}=0,\cr
z_{11}^{(4)}+z_{22}^{(4)}+z_{33}^{(4)}+z_{44}^{(4)}=0,\cr
z_{11}^{(1)}+z_{12}^{(4)}+z_{23}^{(4)}+z_{34}^{(4)}+z_{14}^{(4)}=0.
\end{cases}
\]
Thus, $L=3$ and
\begin{align*}
&C_1^{(1)}=[1],& &C_1^{(2)}=\left[\begin{matrix} 1&0\cr 0&1\end{matrix}\right],& &C_1^{(4)}=0,\\
&C_2^{(1)}=0,& &C_2^{(2)}=0,& &C_2^{(4)}=\left[
\begin{matrix}
1&0&0&0\cr
0&1&0&0\cr
0&0&1&0\cr
0&0&0&1
\end{matrix}\right],\\
&C_3^{(1)}=[1],& &C_3^{(2)}=0,& &C_2^{(4)}=\left[
\begin{matrix}
0&1&0&1\cr
0&0&1&0\cr
0&0&0&1\cr
0&0&0&0
\end{matrix}\right].
\end{align*}
We have
\[
\text{type}(a_1C_1^{(1)}+a_2C_2^{(1)}+a_3C_3^{(1)})=
\begin{cases}
(1,0,0,0)&\text{if}\ a_1=a_3,\cr
(1,0,0,1)&\text{if}\ a_1\ne a_3,
\end{cases}
\]
\[
\text{type}(a_1C_1^{(2)}+a_2C_2^{(2)}+a_3C_3^{(2)})=
\begin{cases}
(2,0,0,0)&\text{if}\ a_1=0,\cr
(2,0,0,1)&\text{if}\ a_1=1,
\end{cases}
\]
\[
\text{type}(a_1C_1^{(4)}+a_2C_2^{(4)}+a_3C_3^{(4)})=
\begin{cases}
(4,0,0,0)&\text{if}\ (a_2,a_3)=(0,0),\cr
(4,1,0,0)&\text{if}\ (a_2,a_3)=(0,1),\cr
(4,0,0,1)&\text{if}\ (a_2,a_3)=(1,0),\cr
(4,1,1,0)&\text{if}\ (a_2,a_3)=(1,1).
\end{cases}
\]
As $(a_1,a_2,a_3)$ runs over $\Bbb F_2^3$, 
\[
\text{type}(\frak C_{(a_1,a_2,a_3)})=
\begin{cases}
(\theta, 0,0,0)& 1\ \text{time},\cr
(\theta, 0,0,\delta(\alpha_4))& 1\ \text{time},\cr
(\theta, \alpha_4,0,\delta(\alpha_2))& 1\ \text{time},\cr
(\theta, \alpha_4,\alpha_4,\delta(\alpha_2))& 1\ \text{time},\cr
(\theta, \alpha_4,0,\delta(\alpha_1))& 1\ \text{time},\cr
(\theta, \alpha_4,\alpha_4,\delta(\alpha_1))& 1\ \text{time},\cr
(\theta, 0,0,\delta(\alpha_1+\alpha_2))& 1\ \text{time},\cr
(\theta, 0,0,\delta(\alpha_1+\alpha_2+\alpha_4))& 1\ \text{time}.
\end{cases}
\]
By \eqref{n(A)},
\[
\begin{split}
\frak n(A_{10})=\,& 2^{\alpha_1+2\alpha_2+4\alpha_4-3}\bigl[4-\delta(\alpha_4)
-\delta(\alpha_1+\alpha_2)-\delta(\alpha_1+\alpha_2+\alpha_4)\bigr]\cr
&\! +2^{\alpha_1+2\alpha_2+3\alpha_4-3}\bigl(1+(-1)^{\alpha_4}\bigr)\bigl(2-\delta(\alpha_1)
- \delta(\alpha_2)\bigr).
\end{split}
\]
By \eqref{2.9} and \eqref{5.2.4},
\[
\begin{split}
&|\text{Fix}(A_{10},P_\lambda)|\cr
=\,& 2^{4\sum_{\nu(i)>2}\lambda_i}\, \frak n(A_{10})\cr
=\,& 2^{4\lambda_{0,8}}\Bigl[2^{\lambda_{1,2}+2\lambda_{2,4}+4\lambda_{4,8}-3}\bigl[4-\delta(\lambda_{4,8})
-\delta(\lambda_{1,2}+\lambda_{2,4})-\delta(\lambda_{1,2}+\lambda_{2,4}+\lambda_{4,8})\bigr]\cr
&\! +2^{\lambda_{1,2}+2\lambda_{2,4}+3\lambda_{4,8}-3}\bigl(1+(-1)^{\lambda_{4,8}}\bigr)\bigl(2-\delta(\lambda_{1,2})
- \delta(\lambda_{2,4})\bigr)\Bigr].
\end{split}
\]


\subsection{Computation of $|{\rm Fix}(A_{11},P_\lambda)|$}\

We have $o(A_{11})=6$, $s_1=s_2=0$, $s_3=2$, $s_6=4$, 
\[
B_1=B_2=\emptyset,\qquad 
B_3=\left[\begin{matrix}
1&0\cr
1&1\cr
0&1\cr
1&0
\end{matrix}\right],\qquad
B_6=I_4
\]
and
\begin{equation}\label{5.10}
\alpha_3=\lambda_{3,6},\qquad \alpha_6=\lambda_{6,12}\qquad \theta=2\alpha_3+4\alpha_6.
\end{equation}
With $Y_dY_d^T=[z_{ij}^{(d)}]$, equation \eqref{2.10} becomes
\[
\begin{split}
&\left[
\begin{matrix}
0&z_{11}^{(3)}+z_{12}^{(3)}+z_{22}^{(3)}& z_{11}^{(3)}+z_{12}^{(3)}+z_{22}^{(3)} & 0\cr
z_{11}^{(3)}+z_{12}^{(3)}+z_{22}^{(3)} & 0 & z_{11}^{(3)}+z_{12}^{(3)}+z_{22}^{(3)} & z_{11}^{(3)}+z_{12}^{(3)}+z_{22}^{(3)}\cr
z_{11}^{(3)}+z_{12}^{(3)}+z_{22}^{(3)} & z_{11}^{(3)}+z_{12}^{(3)}+z_{22}^{(3)} & 0 & z_{11}^{(3)}+z_{12}^{(3)}+z_{22}^{(3)}\cr
0 & z_{11}^{(3)}+z_{12}^{(3)}+z_{22}^{(3)} & z_{11}^{(3)}+z_{12}^{(3)}+z_{22}^{(3)} 0
\end{matrix}\right]+\cr
\cr
&\left[
\begin{matrix}
0 & z_{12}^{(6)}+z_{14}^{(6)}+z_{34}^{(6)}\cr
z_{12}^{(6)}+z_{14}^{(6)}+z_{34}^{(6)} & 0\cr
z_{11}^{(6)}+z_{13}^{(6)}+z_{22}^{(6)}+z_{24}^{(6)}+z_{33}^{(6)}+z_{44}^{(6)}&
z_{12}^{(6)}+z_{14}^{(6)}+z_{34}^{(6)}\cr
0 & z_{11}^{(6)}+z_{13}^{(6)}+z_{22}^{(6)}+z_{24}^{(6)}+z_{33}^{(6)}+z_{44}^{(6)}
\end{matrix}\right.\cr
\cr
&\kern2mm
\left.
\begin{matrix}
z_{11}^{(6)}+z_{13}^{(6)}+z_{22}^{(6)}+z_{24}^{(6)}+z_{33}^{(6)}+z_{44}^{(6)} & 0\cr
z_{12}^{(6)}+z_{14}^{(6)}+z_{34}^{(6)} &
z_{11}^{(6)}+z_{13}^{(6)}+z_{22}^{(6)}+z_{24}^{(6)}+z_{33}^{(6)}+z_{44}^{(6)}\cr
0 & z_{12}^{(6)}+z_{14}^{(6)}+z_{34}^{(6)}\cr
z_{12}^{(6)}+z_{14}^{(6)}+z_{34}^{(6)} & 0
\end{matrix}\right]=0,
\end{split}
\]
i.e.,
\[
\begin{cases}
z_{11}^{(3)}+z_{12}^{(3)}+z_{22}^{(3)}+z_{12}^{(6)}+z_{14}^{(6)}+z_{34}^{(6)}=0,\cr
z_{11}^{(3)}+z_{12}^{(3)}+z_{22}^{(3)}+z_{11}^{(6)}+z_{13}^{(6)}+z_{22}^{(6)}+z_{24}^{(6)}+z_{33}^{(6)}+z_{44}^{(6)}=0.
\end{cases}
\]
So, $L=2$ and
\begin{align*}
&C_1^{(3)}=\left[\begin{matrix} 1&1\cr 0&1 \end{matrix}\right],& &C_1^{(6)}=\left[
\begin{matrix}
0&1&0&1\cr
0&0&0&0\cr
0&0&0&1\cr
0&0&0&0
\end{matrix}\right],\\
&C_2^{(3)}=\left[\begin{matrix} 1&1\cr 0&1 \end{matrix}\right],& &  
C_2^{(6)}=\left[
\begin{matrix}
1&0&1&0\cr
0&1&0&1\cr
0&0&1&0\cr
0&0&0&1
\end{matrix}\right].
\end{align*}
We have
\[
\text{type}(a_1C_1^{(3)}+a_2C_2^{(3)})=
\begin{cases}
(1,0,0,0)&\text{if}\ a_1=a_2,\cr
(1,1,1,0)&\text{if}\ a_1\ne a_2,
\end{cases}
\]
\[
\text{type}(a_1C_1^{(6)}+a_2C_2^{(6)})=
\begin{cases}
(4,0,0,0)&\text{if}\ (a_1,a_2)=(0,0),\cr
(4,2,0,0)&\text{if}\ a_1\ne a_2,\cr
(4,1,1,0)&\text{if}\ (a_1,a_2)=(1,1),
\end{cases}
\]
\[
\text{type}(\frak C_{(a_1,a_2)})=
\begin{cases}
(\theta,0,0,0)&\text{if}\ (a_1,a_2)=(0,0),\cr
(\theta,\alpha_3+2\alpha_6,\alpha_3,0)&\text{if}\ a_1\ne a_2,\cr
(\theta,\alpha_6,\alpha_6,0)&\text{if}\ (a_1,a_2)=(1,1).
\end{cases}
\]
By \eqref{n(A)},
\[
\frak n(A_{11})=2^{2\alpha_3+4\alpha_6-2}+2^{2\alpha_3+3\alpha_6-2}(-1)^{\alpha_6} + 2^{\alpha_3+2\alpha_6-1}(-1)^{\alpha_3}.
\]
By \eqref{2.9} and \eqref{5.10},
\[
\begin{split}
\text{Fix}(A_{11},P_\lambda)=\,& 2^{4\sum_{i\equiv 0\,(12)}\lambda_i}\,\frak n(A_{11})\cr
=\,& 2^{4\lambda_{0,12}-1}\bigl[
2^{2\lambda_{3,6}+4\lambda_{6,12}-1}\cr
&+2^{2\lambda_{3,6}+3\lambda_{6,12}-1}(-1)^{\lambda_{6,12}} + 2^{\lambda_{3,6}+2\lambda_{6,12}}(-1)^{\lambda_{3,6}}\bigr].
\end{split}
\]


\subsection{Computation of $|{\rm Fix}(A_{12},P_\lambda)|$}\
 
We have $o(A_{12})=5$,  $s_1=0$, $s_5=4$,
$B_1=\emptyset$, $B_5=I_4$ and
\begin{equation}\label{5.24}
\alpha_5=\lambda_{5,10}.
\end{equation}
With $Y_dY_d^T=[z_{ij}^{(d)}]$, equation \eqref{2.10} becomes
\[
\begin{split}
&\left[ 
\begin{matrix}
0& z_{11}^{(5)}+z_{13}^{(5)}+z_{23}^{(5)}+z_{24}^{(5)}+z_{44}^{(5)}\cr 
z_{11}^{(5)}+z_{13}^{(5)}+z_{23}^{(5)}+z_{24}^{(5)}+z_{44}^{(5)} &0\cr 
z_{12}^{(5)}+z_{14}^{(5)}+z_{22}^{(5)}+z_{33}^{(5)}+z_{34}^{(5)}& z_{11}^{(5)}+z_{13}^{(5)}+z_{23}^{(5)}+z_{24}^{(5)}+z_{44}^{(5)}\cr 
z_{12}^{(5)}+z_{14}^{(5)}+z_{22}^{(5)}+z_{33}^{(5)}+z_{34}^{(5)}& z_{12}^{(5)}+z_{14}^{(5)}+z_{22}^{(5)}+z_{33}^{(5)}+z_{34}^{(5)}\cr 
\end{matrix}\right. \cr
\cr
&\kern2mm \left.
\begin{matrix}
z_{12}^{(5)}+z_{14}^{(5)}+z_{22}^{(5)}+z_{33}^{(5)}+z_{34}^{(5)}& z_{12}^{(5)}+z_{14}^{(5)}+z_{22}^{(5)}+z_{33}^{(5)}+z_{34}^{(5)}\cr
z_{11}^{(5)}+z_{13}^{(5)}+z_{23}^{(5)}+z_{24}^{(5)}+z_{44}^{(5)} &z_{12}^{(5)}+z_{14}^{(5)}+z_{22}^{(5)}+z_{33}^{(5)}+z_{34}^{(5)}\cr
0 &z_{11}^{(5)}+z_{13}^{(5)}+z_{23}^{(5)}+z_{24}^{(5)}+z_{44}^{(5)} \cr
z_{11}^{(5)}+z_{13}^{(5)}+z_{23}^{(5)}+z_{24}^{(5)}+z_{44}^{(5)}&0\cr
\end{matrix}
\right]=0,
\end{split}
\]
i.e.,
\[
\begin{cases}
z_{11}^{(5)}+z_{13}^{(5)}+z_{23}^{(5)}+z_{24}^{(5)}+z_{44}^{(5)}=0,\cr
z_{12}^{(5)}+z_{14}^{(5)}+z_{22}^{(5)}+z_{33}^{(5)}+z_{34}^{(5)}=0.
\end{cases}
\]
Thus, $L=2$ and
\[
C_1^{(5)}=\left[
\begin{matrix}
1&0&1&0\cr
0&0&1&1\cr
0&0&0&0\cr
0&0&0&1
\end{matrix}\right],\qquad
C_2^{(5)}=\left[
\begin{matrix}
0&1&0&1\cr
0&1&0&0\cr
0&0&1&1\cr
0&0&0&0
\end{matrix}\right].
\]
We have
\[
\text{type}(a_1C_1^{(5)}+a_2C_2^{(5)})=
\begin{cases}
(4,0,0,0)&\text{if}\ (a_1,a_2)=(0,0),\cr
(4,2,1,0)&\text{if}\ (a_1,a_2)\ne(0,0),
\end{cases}
\]
\[
\text{type}(\frak C_{(a_1,a_2)})=
\begin{cases}
(4\alpha_5,0,0,0)&\text{if}\ (a_1,a_2)=(0,0),\cr
(4\alpha_5,2\alpha_5,\alpha_5,0)&\text{if}\ (a_1,a_2)\ne(0,0).
\end{cases}
\]
By \eqref{n(A)}, \eqref{2.9} and \eqref{5.24},
\[
\frak n(A_{12})=2^{4\alpha_5-2}+2^{2\alpha_5-2}\, 3(-1)^{\alpha_5},
\]
\[
\begin{split}
|\text{Fix}(A_{12},P_\lambda)|&=2^{4\sum_{i\equiv 0\,(10)}\lambda_i}\,\frak n(A_{12})\cr
&=2^{4\lambda_{0,10}-2}\bigl[
2^{4\lambda_{5,10}}+2^{2\lambda_{5,10}}\, 3(-1)^{\lambda_{5,10}}\bigr].
\end{split}
\]


\subsection{Computation of $|{\rm Fix}(A_{13},P_\lambda)|$ and $|{\rm Fix}(A_{14},P_\lambda)|$}\

The smallest positive integer $t$ such that $x^4+x+1\mid x^t-1$ is $15$. It follows immediately from Lemma~\ref{L2.10}
that
\begin{equation}\label{5.29}
|{\rm Fix}(A_{13},P_\lambda)|=2^{4\lambda_{0,15}}.
\end{equation}
In the same way,
\begin{equation}\label{5.30}
|{\rm Fix}(A_{14},P_\lambda)|=2^{4\lambda_{0,15}}.
\end{equation}


\section{Computation of $\Psi_{\le 5, n}$}

The group $\text{GL}(5,\Bbb F_2)$ has 27 conjugacy classes. A set of representatives, $A_1,\dots$, $A_{27}$, of the conjugacy classes is given in Table~\ref{Tb3}. In this section,
we give explicit formulas for $|\text{Fix}(A_i,P_\lambda)|$, $1\le i\le 27$. These formulas
are computed by the same method in Sections 4 and 5. In many cases, Lemma~\ref{2.9} and the results of Sections 4 and 5 are used. The details of the computations are
omitted.  

\begin{table}
\caption{Information about $\mathcal C({\rm GL}(5,\Bbb F_2))$}\label{Tb3}
\vskip-5mm
\[
\begin{tabular}{l|l|c}
\hline
representative & elementary divisors& $|{\rm cent}_{{\rm GL}(5,\Bbb F_2)}(\ )|$\\ \hline
$A_1=I_5$ & $x+1,\ x+1,\ x+1,\ x+1,\ x+1$ & $2^{10}\cdot 3^2\cdot 5\cdot 7\cdot 31$ \\ \hline
$A_2=I_3\oplus\left[
\begin{matrix}
0&1\cr
1&0\cr
\end{matrix}\right]$ & $x+1,\ x+1,\ x+1,\ (x+1)^2$ & $2^{10}\cdot 3\cdot 7$ \\ \hline
$A_3=I_2\oplus\left[
\begin{matrix}
0&1&0\cr
0&0&1\cr
1&1&1\cr
\end{matrix}\right]$ & $x+1,\ x+1,\ (x+1)^3$ & $2^7\cdot 3$ \\ \hline
$A_4=I_3\oplus\left[
\begin{matrix}
0&1\cr
1&1\cr
\end{matrix}\right]$ & $x+1,\ x+1,\ x+1,\ x^2+x+1$ & $2^3\cdot 3^2\cdot 7$ \\ \hline
$A_5=I_2\oplus\left[
\begin{matrix}
0&1&0\cr
0&0&1\cr
1&1&0\cr
\end{matrix}\right]$ & $x+1,\ x+1,\ x^3+x+1$ & $2\cdot 3\cdot 7$ \\ \hline
$A_6=I_2\oplus\left[
\begin{matrix}
0&1&0\cr
0&0&1\cr
1&0&1\cr
\end{matrix}\right]$ & $x+1,\ x+1,\ x^3+x^2+1$ & $2\cdot 3\cdot 7$ \\ \hline
$A_7=[1]\oplus\left[
\begin{matrix}
0&1\cr
1&0\cr
\end{matrix}\right]\oplus\left[
\begin{matrix}
0&1\cr
1&0\cr
\end{matrix}\right]$ & $x+1,\ (x+1)^2,\ (x+1)^2$ & $2^9\cdot 3$ \\ \hline
$A_8=[1]\oplus\left[
\begin{matrix}
0&1\cr
1&0\cr
\end{matrix}\right]\oplus\left[
\begin{matrix}
0&1\cr
1&1\cr
\end{matrix}\right]$ & $x+1,\ (x+1)^2,\ x^2+x+1$ & $2^3\cdot 3$ \\ \hline
$A_9=[1]\oplus\left[
\begin{matrix}
0&1\cr
1&1\cr
\end{matrix}\right]\oplus\left[
\begin{matrix}
0&1\cr
1&1\cr
\end{matrix}\right]$ & $x+1,\ x^2+x+1,\ x^2+x+1$ & $2^2\cdot 3^2\cdot 5$ \\ \hline
$A_{10}=[1]\oplus\left[
\begin{matrix}
0&1&0&0\cr
0&0&1&0\cr
0&0&0&1\cr
1&0&0&0\cr
\end{matrix}\right]$ & $x+1,\ (x+1)^4$ & $2^5$ \\ \hline
$A_{11}=[1]\oplus\left[
\begin{matrix}
0&1&0&0\cr
0&0&1&0\cr
0&0&0&1\cr
1&0&1&0\cr
\end{matrix}\right]$ & $x+1,\ (x^2+x+1)^2$ & $2^2\cdot 3$ \\ \hline
$A_{12}=[1]\oplus\left[
\begin{matrix}
0&1&0&0\cr
0&0&1&0\cr
0&0&0&1\cr
1&1&1&1\cr
\end{matrix}\right]$ & $x+1,\ x^4+x^3+x^2+x+1$ & $3\cdot 5$ \\ \hline
$A_{13}=[1]\oplus\left[
\begin{matrix}
0&1&0&0\cr
0&0&1&0\cr
0&0&0&1\cr
1&1&0&0\cr
\end{matrix}\right]$ & $x+1,\ x^4+x+1$ & $3\cdot 5$ \\ \hline
$A_{14}=[1]\oplus\left[
\begin{matrix}
0&1&0&0\cr
0&0&1&0\cr
0&0&0&1\cr
1&0&0&1\cr
\end{matrix}\right]$ & $x+1,\ x^4+x^3+1$ & $3\cdot 5$ \\ \hline
$A_{15}=\left[\begin{matrix} 0&1\cr 1&0\end{matrix}\right]\oplus\left[
\begin{matrix}
0&1&0\cr
0&0&1\cr
1&1&1
\end{matrix}\right]$ & $(x+1)^2,\ (x+1)^3$ & $2^7$ \\ \hline
$A_{16}=\left[\begin{matrix} 0&1\cr 1&0\end{matrix}\right]\oplus\left[
\begin{matrix}
0&1&0\cr
0&0&1\cr
1&1&0
\end{matrix}\right]$ & $(x+1)^2,\ x^3+x+1$ & $2\cdot 7$ \\ \hline
\end{tabular}
\]
\end{table}

\addtocounter{table}{-1}

\begin{table}
\caption{Continued}
\vskip-5mm
\[
\begin{tabular}{l|l|c}
\hline
$A_{17}=\left[\begin{matrix} 0&1\cr 1&0\end{matrix}\right]\oplus\left[
\begin{matrix}
0&1&0\cr
0&0&1\cr
1&0&1
\end{matrix}\right]$ & $(x+1)^2,\ x^3+x^2+1$ & $2\cdot 7$ \\ \hline
$A_{18}=\left[\begin{matrix} 0&1\cr 1&1\end{matrix}\right]\oplus\left[
\begin{matrix}
0&1&0\cr
0&0&1\cr
1&1&1
\end{matrix}\right]$ & $x^2+x+1,\ (x+1)^3$ & $2^2\cdot 3$ \\ \hline
$A_{19}=\left[\begin{matrix} 0&1\cr 1&1\end{matrix}\right]\oplus\left[
\begin{matrix}
0&1&0\cr
0&0&1\cr
1&1&0
\end{matrix}\right]$ & $x^2+x+1,\ x^3+x+1$ & $3\cdot 7$ \\ \hline
$A_{20}=\left[\begin{matrix} 0&1\cr 1&1\end{matrix}\right]\oplus\left[
\begin{matrix}
0&1&0\cr
0&0&1\cr
1&0&1
\end{matrix}\right]$ & $x^2+x+1,\ x^3+x^2+1$ & $3\cdot 7$ \\ \hline
$A_{21}=\left[
\begin{matrix}
0&1&0&0&0\cr
0&0&1&0&0\cr
0&0&0&1&0\cr
0&0&0&0&1\cr
1&1&0&0&1
\end{matrix}\right]$ & $(x+1)^5$ & $2^4$ \\ \hline
$A_{22}=\left[
\begin{matrix}
0&1&0&0&0\cr
0&0&1&0&0\cr
0&0&0&1&0\cr
0&0&0&0&1\cr
1&0&1&0&0
\end{matrix}\right]$ & $x^5+x^2+1$ & $31$ \\ \hline
$A_{23}=\left[
\begin{matrix}
0&1&0&0&0\cr
0&0&1&0&0\cr
0&0&0&1&0\cr
0&0&0&0&1\cr
1&0&0&1&0
\end{matrix}\right]$ & $x^5+x^3+1$ & $31$ \\ \hline
$A_{24}=\left[
\begin{matrix}
0&1&0&0&0\cr
0&0&1&0&0\cr
0&0&0&1&0\cr
0&0&0&0&1\cr
1&1&1&1&0
\end{matrix}\right]$ & $x^5+x^3+x^2+x+1$ & $31$ \\ \hline
$A_{25}=\left[
\begin{matrix}
0&1&0&0&0\cr
0&0&1&0&0\cr
0&0&0&1&0\cr
0&0&0&0&1\cr
1&1&1&0&1
\end{matrix}\right]$ & $x^5+x^4+x^2+x+1$ & $31$ \\ \hline
$A_{26}=\left[
\begin{matrix}
0&1&0&0&0\cr
0&0&1&0&0\cr
0&0&0&1&0\cr
0&0&0&0&1\cr
1&1&0&1&1
\end{matrix}\right]$ & $x^5+x^4+x^3+x+1$ & $31$ \\ \hline
$A_{27}=\left[
\begin{matrix}
0&1&0&0&0\cr
0&0&1&0&0\cr
0&0&0&1&0\cr
0&0&0&0&1\cr
1&0&1&1&1
\end{matrix}\right]$ & $x^5+x^4+x^3+x^2+1$ & $31$ \\ \hline
\end{tabular}
\]
\end{table}

\vskip5mm
\hrule 
\vskip3mm
\[
\begin{split}
|\text{Fix}(A_1,P_\lambda)|=\,& 2^{5\lambda_{0,2}}\Bigl[ 2^{5\lambda_{1,2}-15}\bigl(18260-18259\,\delta(\lambda_{1,2})\bigr)\cr
&+2^{4\lambda_{1,2}-15}\, 155\bigl( 3+(-1)^{\lambda_{1,2}}\bigr)
+2^{3\lambda_{1,2}-12}\, 217\bigl( 5+3(-1)^{\lambda_{1,2}}\bigr)\Bigr].
\end{split}
\]

\vskip3mm
\hrule
\vskip3mm

\[
\begin{split}
|\text{Fix}(A_2,P_\lambda)|=\,& 2^{5\lambda_{0,4}}\Bigl[2^{4\lambda_{1,2}+5\lambda_{2,4}-11}
\bigl[200-99\delta(\lambda_{1,2})-99\delta(\lambda_{1,2}+\lambda_{2,4})-\delta(\lambda_{2,4})\bigr]\cr
&\! +2^{3\lambda_{1,2}+5\lambda_{2,4}-11}\, 7 \bigl(3+(-1)^{\lambda_{1,2}}\bigr)\bigl(2-\delta(\lambda_{2,4})\bigr)\cr
&\! +2^{4\lambda_{1,2}+4\lambda_{2,4}-6}\, 21 \bigl(1-\delta(\lambda_{1,2})\bigr)
+ 2^{3\lambda_{1,2}+4\lambda_{2,4}-8}\, 7 \bigl(3+(-1)^{\lambda_{1,2}}\bigr)\cr
&\! + 2^{2\lambda_{1,2}+4\lambda_{2,4}-7}\, 7 \bigl(5+3(-1)^{\lambda_{1,2}}\bigr)\Bigr].
\end{split}
\]

\vskip3mm
\hrule
\vskip3mm

\[
\begin{split}
|\text{Fix}(A_3,P_\lambda)|=\,& 2^{5\lambda_{0,8}}\Bigl[
2^{3\lambda_{1,2}+4\lambda_{2,4}+5\lambda_{4,8}-7}
\bigl[24-3\delta(\lambda_{1,2})
-8\delta(\lambda_{1,2}+\lambda_{2,4})\cr
&\! -3\delta(\lambda_{1,2}+\lambda_{2,4}+\lambda_{4,8})-8\delta(\lambda_{1,2}+\lambda_{4,8})-\delta(\lambda_{2,4}+\lambda_{4,8})\bigr]\cr
&\! +2^{3\lambda_{1,2}+3\lambda_{2,4}+5\lambda_{4,8}-4} \, 3 \bigl[2-\delta(\lambda_{1,2})-\delta(\lambda_{1,2}+\lambda_{4,8})\bigr]\cr
&\! +2^{2\lambda_{1,2}+4\lambda_{2,4}+5\lambda_{4,8}-7}\bigl( 3+(-1)^{\lambda_{1,2}}\bigr) \bigl(2- \delta(\lambda_{2,4}+\lambda_{4,8})\bigr)\cr
&\! +2^{2\lambda_{1,2}+3\lambda_{2,4}+5\lambda_{4,8}-6}\, 3 \bigl( 3+(-1)^{\lambda_{1,2}}\bigr) \bigl(2- \delta(\lambda_{4,8})\bigr)\Bigr].
\end{split}
\]

\vskip3mm
\hrule
\vskip3mm

\[
\begin{split}
|\text{Fix}(A_4,P_\lambda)|=\,& 2^{3\lambda_{0,2}+2\lambda_{0,6}-1}\bigl(2^{2\lambda_{3,6}}
+(-1)^{\lambda_{3,6}}2^{\lambda_{3,6}}\bigr)\cr
&\cdot \bigr[2^{3\lambda_{1,2}-6}\bigl(36-35\, \delta(\lambda_{1,2})\bigr)+2^{2\lambda_{1,2}-6}\, 7
\bigl( 3+(-1)^{\lambda_{1,2}}\bigr)\bigr].
\end{split}
\]

\vskip3mm
\hrule
\vskip3mm

For $j=5,6$,
\[
|\text{Fix}(A_j,P_\lambda)|=2^{2\lambda_{0,2}+3\lambda_{0,7}}\bigl[
2^{2\lambda_{1,2}-3}\bigl(4-3\,\delta(\lambda_{1,2})\bigr)+2^{\lambda_{1,2}-3}\bigl(3+(-1)^{\lambda_{1,2}}\bigr)\bigr].
\]

\vskip3mm
\hrule
\vskip3mm

\[
\begin{split}
|\text{Fix}(A_7,P_\lambda)|=\,& 2^{5\lambda_{0,4}}\Bigl[
2^{3\lambda_{1,2}+5\lambda_{2,4}-9}\bigl[80-7\delta(\lambda_{1,2})-69\delta(\lambda_{1,2}+
\lambda_{2,4})-3\delta(\lambda_{2,4})\bigr]\cr
& +2^{3\lambda_{1,2}+4\lambda_{2,4}-9}\bigl[72 +8(-1)^{\lambda_{2,4}}-7(-1)^{\lambda_{2,4}}
\delta(\lambda_{1,2})-69\delta(\lambda_{1,2})\bigr]\cr
& +2^{2\lambda_{1,2}+5\lambda_{2,4}-7}\, 3 \bigl(3+(-1)^{\lambda_{1,2}} \bigr)
\bigl(1- \delta(\lambda_{2,4})\bigr)\cr
&+2^{2\lambda_{1,2}+4\lambda_{2,4}-7}\, 3 \bigl(3+(-1)^{\lambda_{1,2}} \bigr)
+2^{3\lambda_{1,2}+3\lambda_{2,4}-2}\bigl(1- \delta(\lambda_{1,2})\bigr)\cr
&+2^{2\lambda_{1,2}+3\lambda_{2,4}-4}\bigl(3+(-1)^{\lambda_{1,2}} \bigr)\Bigr]. 
\end{split}
\]

\vskip3mm
\hrule
\vskip3mm

\[
\begin{split}
|\text{Fix}(A_8,P_\lambda)|=\,& 2^{3\lambda_{0,4}+2\lambda_{0,6}-1}\bigl(2^{2\lambda_{3,6}}
+(-1)^{\lambda_{3,6}}2^{\lambda_{3,6}}\bigr)\cr
&\cdot\Bigl[2^{2\lambda_{1,2}+3\lambda_{2,4}-4}\bigl[8-3\delta(\lambda_{1,2})-
\delta(\lambda_{2,4})-3\delta(\lambda_{1,2}+\lambda_{2,4})\bigr]\cr
&+2^{\lambda_{1,2}+2\lambda_{2,4}-3}\bigl(3+ (-1)^{\lambda_{1,2}}\bigr)\Bigr].
\end{split}
\]

\vskip3mm
\hrule
\vskip3mm

\[
\begin{split}
|\text{Fix}(A_9,P_\lambda)|=\,& 2^{\lambda_{0,1}+4\lambda_{0,6}-4}\bigl(2-\delta(\lambda_{1,2}) \bigr)\cr
&\cdot
\bigl[2^{4\lambda_{3,6}-1}+2^{3\lambda_{3,6}-1}\, 5 (-1)^{\lambda_{3,6}}+
2^{2\lambda_{3,6}}\, 5\bigr].
\end{split}
\]

\vskip3mm
\hrule
\vskip3mm
 
\[
\begin{split}
|\text{Fix}(A_{10},P_\lambda)|=\,& 2^{5\lambda_{0,8}}\Bigl[2^{2\lambda_{1,2}+3\lambda_{2,4}+
5\lambda_{4,8}-5}\bigl[8-\delta(\lambda_{1,2})-2\delta(\lambda_{1,2}+\lambda_{2,4})\cr
& -2\delta(\lambda_{1,2}+\lambda_{2,4}+\lambda_{4,8})-\delta(\lambda_{1,2}+\lambda_{4,8})
-\delta(\lambda_{4,8})\bigr]\cr
&+2^{2\lambda_{1,2}+3\lambda_{2,4}+4\lambda_{4,8}-5}\bigl(1+(-1)^{\lambda_{4,8}} \bigr)\cr
&\cdot
\bigl[ 4-2\delta(\lambda_{1,2})-\delta(\lambda_{1,2}+\lambda_{2,4})-\delta(\lambda_{2,4})
\bigr]\cr
&+2^{\lambda_{1,2}+2\lambda_{2,4}+5\lambda_{4,8}-4}\bigl(3+ (-1)^{\lambda_{1,2}} \bigr)
\bigl(1-\delta(\lambda_{4,8}) \bigr)\cr
&+2^{\lambda_{1,2}+2\lambda_{2,4}+4\lambda_{4,8}-4}\bigl(3+ (-1)^{\lambda_{1,2}} \bigr)\Bigr].
\end{split}
\]
 
\vskip3mm
\hrule
\vskip3mm
 
\[
\begin{split}
|\text{Fix}(A_{11},P_\lambda)|=\,& 2^{\lambda_{0,1}+4\lambda_{0,12}-2}\bigl(2-\delta(\lambda_{1,2})\bigr)
\bigl[2^{2\lambda_{3,6}+4\lambda_{6,12}-1}\cr
&+2^{2\lambda_{3,6}+3\lambda_{6,12}-1}(-1)^{\lambda_{6,12}}
+2^{\lambda_{3,6}+2\lambda_{6,12}}(-1)^{\lambda_{3,6}}\bigr].
\end{split}
\]
 
\vskip3mm
\hrule
\vskip3mm
 
\[
|\text{Fix}(A_{12},P_\lambda)|=2^{\lambda_{0,1}+4\lambda_{0,10}-3}\bigl(2-\delta(\lambda_{1,2})\bigr)\bigl[2^{4\lambda_{5,10}}+2^{2\lambda_{5,10}}\, 3 
(-1)^{\lambda_{5,10}}\bigr].
\]
 
\vskip3mm
\hrule
\vskip3mm

For $j=13,14$,
\[
|\text{Fix}(A_j,P_\lambda)|=2^{\lambda_{0,1}+4\lambda_{0,15}-1}\bigl(2-\delta(\lambda_{1,2})\bigr).
\]
 
\vskip3mm
\hrule
\vskip3mm

\[
\begin{split}
|\text{Fix}(A_{15},P_\lambda)|=\,& 2^{5\lambda_{0,8}}\Bigl[2^{2\lambda_{1,2}+4\lambda_{2,4}+
5\lambda_{4,8}-6}\bigl[16-\delta(\lambda_{1,2})-5\delta(\lambda_{1,2}+\lambda_{2,4})\cr
& -4\delta(\lambda_{1,2}+\lambda_{2,4}+\lambda_{4,8})-2\delta(\lambda_{1,2}+\lambda_{4,8})\cr
&-\delta(\lambda_{2,4})
-2\delta(\lambda_{2,4}+\lambda_{4,8})\bigr]\cr
& +2^{2\lambda_{1,2}+3\lambda_{2,4}+5\lambda_{4,8}-6}\bigl[12+4(-1)^{\lambda_{2,4}}-2(-1)^{\lambda_{2,4}}\delta(\lambda_{1,2})\cr
&-4\delta(\lambda_{1,2})
-5\delta(\lambda_{1,2}+\lambda_{4,8})-(-1)^{\lambda_{2,4}}\delta(\lambda_{1,2}+\lambda_{4,8})\cr
&-\delta(\lambda_{4,8})-(-1)^{\lambda_{2,4}}\delta(\lambda_{4,8})\bigr]\cr
& +2^{\lambda_{1,2}+2\lambda_{2,4}+4\lambda_{4,8}-3}\bigl(3+(-1)^{\lambda_{1,2}}\bigr)
\Bigr].
\end{split}
\]

\vskip3mm
\hrule
\vskip3mm

For $j=16,17$,
\[
|\text{Fix}(A_j,P_\lambda)|=2^{\lambda_{1,2}+2\lambda_{0,2}+3\lambda_{0,7}-2}
\bigl[ 4-\delta(\lambda_{1,2})-\delta(\lambda_{2,4})-\delta(\lambda_{1,2}+\lambda_{2,4})\bigr].
\]
 
\vskip3mm
\hrule
\vskip3mm

\[
\begin{split}
|\text{Fix}(A_{18},P_\lambda)|=\,& 2^{2\lambda_{0,6}+\lambda_{1,2}+2\lambda_{2,4}
+3\lambda_{0,4}-3}\bigl[ 2^{2\lambda_{3,6}}+(-1)^{\lambda_{3,6}}2^{\lambda_{3,6}}\bigr]\cr
&\cdot\bigl[ 4-\delta(\lambda_{1,2}+\lambda_{2,4})-\delta(\lambda_{1,2}+\lambda_{4,8})
-\delta(\lambda_{2,4}+\lambda_{4,8})\bigr].
\end{split}
\]
 
\vskip3mm
\hrule
\vskip3mm

For $j=19,20$,
\[
|\text{Fix}(A_j,P_\lambda)|=2^{2\lambda_{0,6}+3\lambda_{0,7}-1}\bigl[2^{2\lambda_{3,6}}
+(-1)^{\lambda_{3,6}}2^{\lambda_{3,6}}\bigr].
\]
 
\vskip3mm
\hrule
\vskip3mm

\[
\begin{split}
|\text{Fix}(A_{21},P_\lambda)|=\,& 2^{5\lambda_{0,16}}\Bigl[
2^{\lambda_{1,2}+2\lambda_{2,4}+4\lambda_{4,8}+5\lambda_{8,16}-3}\bigl[ 
4-\delta(\lambda_{1,2}+\lambda_{2,4})\cr
& -\delta(\lambda_{1,2}+\lambda_{2,4}+\lambda_{4,8})
-\delta(\lambda_{4,8})\bigr]\cr
& +2^{\lambda_{1,2}+2\lambda_{2,4}+3\lambda_{4,8}+5\lambda_{8,16}-3}
\bigl(1+(-1)^{\lambda_{4,8}}\bigr)\cr
&\cdot \bigl[2-\delta(\lambda_{1,2}+\lambda_{8,16})-\delta(\lambda_{2,4}+\lambda_{8,16})
\bigr]\Bigr].
\end{split}
\]
 
\vskip3mm
\hrule
\vskip3mm

For $22\le j\le 27$,
\[
|\text{Fix}(A_j,P_\lambda)|=2^{5\lambda_{0,31}}.
\]
 
\vskip3mm
\hrule
\vskip3mm


\section{Numerical Results}

Two tables are included in this section. Table~\ref{Tb4} contains the values of $\Psi_{\le k, n}$ with $k\le 5$ and $n\le 40$. Table~\ref{Tb5}
contains the values of $\Psi_{k,n}$ with $k$ and $n$ in the same ranges

\begin{table}
\caption{Values of $\Psi_{\le k, n}$, $k\le 5$, $n\le 40$}\label{Tb4}
\vskip-5mm
\[
\begin{tabular}{c||r|r|r|r|r|r}
\hline
$n\backslash k$ & 0 & \multicolumn{1}{c|} 1 & \multicolumn{1}{c|} 2 & \multicolumn{1}{c|}3 
& \multicolumn{1}{c|} 4  & \multicolumn{1}{c} 5 \\ \hline \hline
1 & 1 & 1 & 1 & 1 & 1 & 1  \\ \hline
2 & 1 & 2 & 2 & 2 & 2 & 2  \\ \hline
3 & 1 & 2 & 2 & 2 & 2 & 2  \\ \hline
4 & 1 & 3 & 4 & 4 & 4 & 4  \\ \hline
5 & 1 & 3 & 4 & 4 & 4 & 4  \\ \hline
6 & 1 & 4 & 7 & 8 & 8 & 8  \\ \hline
7 & 1 & 4 & 7 & 9 & 9 & 9  \\ \hline
8 & 1 & 5 & 11 & 16 & 18 & 18  \\ \hline
9 & 1 & 5 & 11 & 17 & 20 & 20  \\ \hline
10 & 1 & 6 & 16 & 28 & 37 & 39  \\ \hline
11 & 1 & 6 & 16 & 30 & 42 & 46  \\ \hline
12 & 1 & 7 & 23 & 49 & 77 & 92  \\ \hline
13 & 1 & 7 & 23 & 53 & 89 & 112  \\ \hline
14 & 1 & 8 & 31 & 82 & 157 & 218  \\ \hline
15 & 1 & 8 & 31 & 89 & 187 & 281  \\ \hline
16 & 1 & 9 & 41 & 133 & 323 & 551  \\ \hline
17 & 1 & 9 & 41 & 144 & 389 & 740  \\ \hline
18 & 1 & 10 & 53 & 210 & 654 & 1,447  \\ \hline
19 & 1 & 10 & 53 & 229 & 804 & 2,059  \\ \hline
20 & 1 & 11 & 67 & 325 & 1,324 & 4,029  \\ \hline
21 & 1 & 11 & 67 & 354 & 1,651 & 6,032  \\ \hline
22 & 1 & 12 & 83 & 490 & 2,654 & 11,774  \\ \hline
23 & 1 & 12 & 83 & 534 & 3,356 & 18,581  \\ \hline
24 & 1 & 13 & 102 & 727 & 5,291 & 36,239  \\ \hline
25 & 1 & 13 & 102 & 793 & 6,759 & 59,798  \\ \hline
26 & 1 & 14 & 123 & 1,058 & 10,433 & 116,020  \\ \hline
27 & 1 & 14 & 123 & 1,154 & 13,444 & 198,489  \\ \hline
28 & 1 & 15 & 147 & 1,515 & 20,363 & 382,272  \\ \hline
29 & 1 & 15 & 147 & 1,651 & 26,384 & 670,031  \\ \hline
30 & 1 & 16 & 174 & 2,136 & 39,229 & 1,276,454  \\ \hline
31 & 1 & 16 & 174 & 2,329 & 51,025 & 2,267,431  \\ \hline
32 & 1 & 17 & 204 & 2,972 & 74,574 & 4,260,828  \\ \hline
33 & 1 & 17 & 204 & 3,237 & 97,143 & 7,596,889  \\ \hline
34 & 1 & 18 & 237 & 4,078 & 139,660 & 14,050,410  \\ \hline
35 & 1 & 18 & 237 & 4,439 & 181,923 & 24,965,555  \\ \hline
36 & 1 & 19 & 274 & 5,532 & 257,592 & 45,384,782  \\ \hline
37 & 1 & 19 & 274 & 6,017 & 335,029 & 79,965,507  \\ \hline
38 & 1 & 20 & 314 & 7,418 & 467,600 & 142,792,476  \\ \hline
39 & 1 & 20 & 314 & 8,061 & 606,613 & 248,697,834  \\ \hline
40 & 1 & 21 & 358 & 9,843 & 835,392 & 497,412,483  \\ \hline
\end{tabular}
\]
\end{table}

\begin{table}
\caption{Values of $\Psi_{k, n}$, $k\le 5$, $n\le 40$}\label{Tb5}
\vskip-5mm
\[
\begin{tabular}{c||r|r|r|r|r|r}
\hline
$n\backslash k$ & 0  & \multicolumn{1}{c|} 1 & \multicolumn{1}{c|} 2 & \multicolumn{1}{c|} 3  
& \multicolumn{1}{c|} 4 &\multicolumn{1}{c}5\\ \hline \hline
1 & 1 & 0 & 0 & 0 & 0 & 0  \\ \hline 
2 & 1 & 1 & 0 & 0 & 0 & 0  \\ \hline 
3 & 1 & 1 & 0 & 0 & 0 & 0  \\ \hline 
4 & 1 & 2 & 1 & 0 & 0 & 0  \\ \hline 
5 & 1 & 2 & 1 & 0 & 0 & 0  \\ \hline 
6 & 1 & 3 & 3 & 1 & 0 & 0  \\ \hline 
7 & 1 & 3 & 3 & 2 & 0 & 0  \\ \hline 
8 & 1 & 4 & 6 & 5 & 2 & 0  \\ \hline 
9 & 1 & 4 & 6 & 6 & 3 & 0  \\ \hline 
10 & 1 & 5 & 10 & 12 & 9 & 2  \\ \hline 
11 & 1 & 5 & 10 & 14 & 12 & 4  \\ \hline 
12 & 1 & 6 & 16 & 26 & 28 & 15  \\ \hline 
13 & 1 & 6 & 16 & 30 & 36 & 23  \\ \hline 
14 & 1 & 7 & 23 & 51 & 75 & 61  \\ \hline 
15 & 1 & 7 & 23 & 58 & 98 & 94  \\ \hline 
16 & 1 & 8 & 32 & 92 & 190 & 228  \\ \hline 
17 & 1 & 8 & 32 & 103 & 245 & 351  \\ \hline 
18 & 1 & 9 & 43 & 157 & 444 & 793  \\ \hline 
19 & 1 & 9 & 43 & 176 & 575 & 1,255  \\ \hline 
20 & 1 & 10 & 56 & 258 & 999  & 2,705  \\ \hline 
21 & 1 & 10 & 56 & 287 & 1,297 & 4,381  \\ \hline 
22 & 1 & 11 & 71 & 407 & 2,164 & 9,120  \\ \hline 
23 & 1 & 11 & 71 & 451 & 2,822 & 15,225  \\ \hline 
24 & 1 & 12 & 89 & 625 & 4,564 & 30,948  \\ \hline 
25 & 1 & 12 & 89 & 691 & 5,966 & 53,039  \\ \hline 
26 & 1 & 13 & 109 & 935 & 9,375 & 105,587  \\ \hline 
27 & 1 & 13 & 109 & 1,031 & 12,290 & 185,045  \\ \hline 
28 & 1 & 14 & 132 & 1,368 & 18,848 & 361,909  \\ \hline 
29 & 1 & 14 & 132 & 1,504 & 24,733 & 643,647  \\ \hline 
30 & 1 & 15 & 158 & 1,962 & 37,093 & 1,237,225  \\ \hline 
31 & 1 & 15 & 158 & 2,155 & 48,696 & 2,216,406  \\ \hline 
32 & 1 & 16 & 187 & 2,768 & 71,602 & 4,186,254  \\ \hline 
33 & 1 & 16 & 187 & 3,033 & 93,906 & 7,499,746  \\ \hline 
34 & 1 & 17 & 219 & 3,841 & 135,582 & 13,910,750  \\ \hline 
35 & 1 & 17 & 219 & 4,202 & 177,484 & 24,783,632  \\ \hline 
36 & 1 & 18 & 255 & 5,258 & 252,060 & 45,127,190  \\ \hline 
37 & 1 & 18 & 255 & 5,743 & 329,012 & 79,630,478  \\ \hline 
38 & 1 & 19 & 294 & 7,104 & 460,182 & 142,324,876  \\ \hline 
39 & 1 & 19 & 294 & 7,747 & 598,552 & 248,091,221  \\ \hline 
40 & 1 & 20 & 337 & 9,485 & 825,549 & 496,577,091  \\ \hline 
\end{tabular}
\]
\end{table}



\begin{thebibliography}{99}

\bibitem{Con80}
J. H. Conway and V. Pless, 
{\it On the enumeration of self-dual codes},
J. Combin. Theory A {\bf 28} (1980), 26 -- 53.

\bibitem{Con92}
J. H. Conway, V. Pless, and N. J. A. Sloane, 
{\it The binary self-dual codes of length up to 32: a revised enumeration},
J. Combin. Theory A {\bf 60} (1992), 183 -- 195.

\bibitem{Dic58}
L. E. Dickson,
{\it Linear Groups}, Dover, New York, 1958.

\bibitem{Dil04}
J. F. Dillon and H. Dobbertin,
{\it New cyclic difference set with Singer parameters},
Finite Fields Appl. {\bf 10} (2004), 342 -- 389.

\bibitem{Hou94} X. Hou, 
{\it Classification of cosets of the Reed-Muller code
$R(m-3,m)$}, Discrete Math. {\bf 128} (1994), 203 -- 224.

\bibitem{Hou96}
X. Hou,
{\it ${\rm GL}(m,2)$ acting on $R(r,m)/R(r-1,m)$}, Discrete Math. {\bf 149} (1996), 99 -- 122.

\bibitem{Hou05}
X. Hou, {\it Affinity of permutations of ${\Bbb F}_2^n$}, Discrete Appl.  Math. {\bf 154} (2006), 313 -- 325.

\bibitem{Hou1}
X. Hou,
{\it On the asymptotic number of inequivalent binary self-dual codes},
J. Combin. Theory A, Available online 22 August 2006.

\bibitem{Mac81}
F. J. MacWilliams and N. J. A. Sloane 
{\it The Theory of Error-Correcting Codes},
North-Holland,
New York, 1981

\bibitem{Ple72}
V. Pless,
{\it A classification of self-orthogonal codes over ${\rm GF}(2)$}, Discrete Math. {\bf 3} (1972), 209 -- 246.

\bibitem{Ple75}
V. Pless and N. J. A. Sloane,
{\it On the classification and enumeration of self-dual codes}, J. Combin. Theory A {\bf 18} (1975), 313 -- 335.

\end{thebibliography}
\end{document}